\documentclass{article}

\usepackage[english]{babel}
\usepackage[utf8]{inputenc}
\usepackage{amsmath}
\usepackage{amssymb}
\usepackage[normalem]{ulem}
\usepackage{amsthm}
\usepackage{comment}
\usepackage{faktor}
\usepackage{cancel}
\usepackage[toc,page]{appendix}
\usepackage{hyperref}
\usepackage{amsfonts}
\usepackage{mathrsfs}
\usepackage{physics}
\usepackage{graphicx}
\usepackage{xcolor}
\usepackage[colorinlistoftodos]{todonotes}
\numberwithin{equation}{section}
\newtheorem{theorem}{Theorem}[section]
\newtheorem{corollary}[theorem]{Corollary}
\newtheorem{lemma}[theorem]{Lemma}
\newtheorem{prop}[theorem]{Proposition}
\newtheorem{deff}[theorem]{Definition}
\newtheorem{remark}[theorem]{Remark}
\theoremstyle{definition}

\newcommand{\ray}{\operatorname{Ray}}
\newcommand{\orr}{\operatorname{Orr}}
\newcommand{\lpl}{\left\langle}
\newcommand{\rer}{\right\rangle}

\renewcommand{\theta}{\vartheta}

\newcommand{\eps}{\varepsilon}
\newcommand{\mc}{\mathcal}

\renewcommand{\l}{\left(}
\renewcommand{\r}{\right)}

\renewcommand{\lll}{\left\{}
\newcommand{\rrr}{\right\}}

\newcommand{\Img}{\operatorname{Img}}

\usepackage{arxiv}

\usepackage[utf8]{inputenc} % allow utf-8 input
\usepackage[T1]{fontenc}    % use 8-bit T1 fonts
\usepackage{hyperref}       % hyperlinks
\usepackage{url}            % simple URL typesetting
\usepackage{booktabs}       % professional-quality tables
\usepackage{amsfonts}       % blackboard math symbols
\usepackage{nicefrac}       % compact symbols for 1/2, etc.
\usepackage{microtype}      % microtypography
\usepackage{lipsum}		% Can be removed after putting your text content
\usepackage{graphicx}
\usepackage{doi}

\title{The adjoint Rayleigh and Orr-Sommerfeld equations: Green function and eigenmodes}
\date{}
\author{{Lorenzo Quarisa}\thanks{LQ is supported by the Warwick Mathematics Institute Centre for Doctoral Training, and gratefully acknowledges funding by University of Warwick’s EU Chancellors' Scholarship scheme. } \\
	Mathematics Institute\\
	University of Warwick\\
	Coventry CV47AL, United Kingdom \\
	\texttt{lorenzo.quarisa@warwick.ac.uk}  \\
	\And
{José L. Rodrigo} \\
	Mathematics Institute\\
	University of Warwick\\
	Coventry CV47AL, United Kingdom \\
	\texttt{j.rodrigo@warwick.ac.uk} \\
}

\renewcommand{\shorttitle}{The adjoint Rayleigh and Orr-Sommerfeld equations: Green function and eigenmodes}

\hypersetup{
pdftitle={he adjoint Rayleigh and Orr-Sommerfeld equations: Green function and eigenmodes},
pdfauthor={David S.~Hippocampus, Elias D.~Striatum},
pdfkeywords={Boundary layers, Prandtl expansion, zero viscosity limit, Navier friction boundary condition},
}

\begin{document}
   
 \today
    
\maketitle

\begin{abstract}
	
The Rayleigh and Orr-Sommerfeld equations are ODEs which arise from the linearized Euler and Navier-Stokes equation around a shear flow. In this paper, we consider the adjoints of the Rayleigh and Orr-Sommerfeld equations on $[0,\infty)$ with respect to the complex $L^2$ product. In the viscous case, we consider a family of viscosity-dependent Navier boundary conditions, which in the limit corresponds to the no-slip condition. We rigorously establish existence and asymptotic properties of their eigenvalues, eigenmodes and Green functions away from critical layers. The adjoint operators are useful because they also allow us to deduce properties about the kernels and images of the original operators.

\end{abstract}

% keywords can be removed

\tableofcontents 
\section{Introduction}
The Rayleigh and Orr-Sommerfeld equations are ordinary differential equations which appear in fluid dynamics, resulting from the linearization of the Euler and Navier-Stokes equations respectively, under the assumption that the solution $u(x,y,t)$ is of the form 
\begin{equation}\label{eq:stream}u(x,y,t):=(\partial_y,-\partial_x)\phi(y)e^{i\alpha(x-ct)}, \end{equation}
for some stream function $\phi$. Under the above ansatz, if we let $U_s:\mathbb{R}_+\to \mathbb{R}$ be a shear flow, $c\in \mathbb{C}$, $\alpha \in \mathbb{Z}\setminus \lll 0\rrr$ then the 2D linearized Euler equation 
 around $U_s$ on $\mathbb{T}\times \mathbb{R}_+$ simplifies to the Rayleigh equation on the half line:
\begin{equation}\label{eq:ray}
\ray_c(\phi):=(U_s-c)\Delta_{\alpha}\phi -U_s''\phi=0,\qquad \Delta_{\alpha}:=\partial_{yy}-\alpha^2,
\end{equation}while for a viscosity $\nu>0$, the linearized Navier-Stokes equations simplify to the Orr-Sommerfeld equation on the half line:
\begin{equation}\label{eq:orrs}
    \orr_{c,\nu}(\phi):=(U_s-c)\Delta_{\alpha}\phi-U_s''\phi-\frac{\nu}{i\alpha}\Delta_{\alpha}^2\phi=0.
\end{equation}

In particular, if $\Im c>0$ then the solution defined by \eqref{eq:stream}  is exponentially growing in time, and we call $U_s$ a \emph{linearly unstable} shear flow. 

These equations were first introduced by Rayleigh \cite{rayleigh}, and Orr and Sommerfeld \cite{orrsomm} respectively. Whilst they are derived from linearized equations, they can also be used to obtain instability results for the corresponding nonlinear equations. In particular, E.Grenier and T.T.Nguyen \cite{greniernguyen} used an Orr-Sommerfeld solution to construct a family of $L^{\infty}$-unstable solutions to the Navier-Stokes equations. Much of this work stems from the effort by the authors to generalize this result to a viscosity dependent Navier boundary condition.

We make some standard assumptions on the shear flow $U_s$. Namely, we required that it is smooth, bounded, vanishes at the origin with its first two derivatives, has a finite limit as $y\to+\infty$ and its derivatives decay exponentially at infinity:
\begin{equation}\label{eq:eta}|U_s^{(k)}(y)|\leq C_ke^{-\eta_0 y},\end{equation}
for some $\eta_0>0$, and for some $C_k>0$ depending on $k \in\mathbb{N}$.
Moreover, we assume that we remain away from any \emph{critical layers}, namely that $c$ and $U_s$ satisfy
\begin{equation}\label{eq:critlay} \inf_{y\geq 0}|U_s-c|>0.\end{equation}
Note that this condition is automatically satisfied if $\Im c\neq 0$. See \cite{critlayer} for an analysis of the linearized Navier-Stokes equations when this hypothesis fails.

Equation \eqref{eq:ray} is naturally coupled with the boundary condition $\phi(0)=0$. In the viscous case \eqref{eq:orrs}, we are interested in the following  Navier boundary condition, for some $\gamma\in \mathbb{R}$:
\begin{equation}\label{eq:orrbc}\begin{cases}
    \phi(0)=0;\\
    \phi'(0)=\nu^{\gamma}\phi''(0).\end{cases}
\end{equation}
This may be seen as a viscosity-dependent Navier friction boundary condition, where we may recover the no-slip case by taking the limit $\gamma \to \infty$.

In this paper, we introduce the adjoints of the above equations, which can be obtained through integration by parts:
\begin{align}\label{eq:raystar}
\ray^*_c(\phi)&=\Delta_{\alpha}(U_s-\bar{c})\phi-U_s''\phi;
\end{align}
and
\begin{equation*}
\orr^*_{c,\nu}(\phi)= \Delta_{\alpha}(U_s-\bar{c})\phi-U_s''\phi+\frac{\nu}{i\alpha}\Delta_{\alpha}^2\phi. 
\end{equation*}

Having fixed an $\alpha \in \mathbb{Z}_+$, 
an \emph{eigenvalue} for these equations is a number $c$ - or a pair $(c,\nu)$ in the viscous case - with $\Im c>0$ and $\nu>0$ such that there exists a solution satisfying the boundary conditions and decay at infinity, which we call an \emph{eigenmode}.

The paper is organized as follows. In Section \ref{sec2} we introduce the adjoint Rayleigh operator $\ray^*_c$ and prove that it shares the same eigenvalues as the original Rayleigh operator. This follows from the fact that if $\phi$ is a solution \eqref{eq:ray}, then
$\bar{\phi}(U_s-\bar{c})^{-1}$ is a solution to the adjoint Rayleigh equation. In Section \ref{sec3} we introduce the adjoint Orr-Sommerfeld equation and construct two pairs of approximate solutions, with different behaviors at infinity. For the slow solutions, we introduce new techniques which allow us to prove their existence and expansions in powers of $\nu$ for all values of $c$, including near Rayleigh eigenvalues. In Section \ref{sec4} we use these approximate solutions to construct a Green function for the adjoint Orr-Sommerfeld equation, which can be made to satisfy the boundary condition provided $(c,\nu)$ is not an eigenvalue. 
Using Rouché's theorem from complex analysis,  we then rigorously prove the existence and asymptotic properties of the eigenvalues of the adjoint Orr-Sommerfeld equation near Rayleigh eigenvalues, and that these eigenvalues in fact match with the eigenvalues of the original equation. Finally, we show as an application that the eigenmodes do not belong to the image of the corresponding Orr-Sommerfeld operator.

The analysis of the fundamental solutions and Green functions is inspired from the work of E. Grenier and T.T.Nguyen     \cite{grenieros} on the standard Rayleigh and Orr-Sommerfeld operators with the no-slip condition.  
Whilst we mainly focus on a viscosity-dependent Navier boundary condition, the result can extend to the no-slip case by taking the limit $\gamma \to \infty$. The main results of the paper can be summarized by the following statement.

\begin{theorem} Fix $\alpha \in \mathbb{Z}_+$.
\begin{enumerate}
\item A pair $(c,\nu)$, with $\Im c>0$ and $\nu>0$ is an adjoint Orr-Sommerfeld eigenvalue if and only if it is an Orr-Sommerfeld eigenvalue.

\item Let $c_0$ be a Rayleigh eigenvalue. Then there exists an integer $\kappa>0$, a constant $C>0$ such that for all $\nu>0$ small enough there exists exactly $\kappa$ (adjoint) Orr-Sommerfeld eigenvalues $(c_{\nu}^1,\nu),\dots,(c_{\nu}^{\kappa},\nu)$ (counted with their multiplicity) with 
$$|c_{\nu}^j-c_0|^{\kappa}\sim C|\mc{O}_{\gamma}(\nu)|\sim C_1\nu^{1/2}\left|\frac{1+C_2\alpha\nu^{\gamma}}{1+C_3\alpha\nu^{\gamma-1/2}}\right|,\qquad \nu \to 0^+,\qquad j=1,\dots,\kappa, $$
for some constants $C_1,C_2,C_3\in\mathbb{C}$ with $C_1\neq 0$ independent from $\nu$ and $\alpha$.
In particular, the corresponding (adjoint) Orr-Sommerfeld eigenmode is of the form
\begin{equation}
    \phi_{c_{\nu}^*,\nu}^*\sim  \frac{\overline{\phi_{c_{\nu},\nu}}}{U_s-\overline{c_{\nu}}}+O\l \mc{O}_{\gamma}(\nu)\r.
\end{equation}
\item Let $(c,\nu)$ be an eigenvalue, and $\phi_{c,\nu}$, $\phi_{c,\nu}^*$ the corresponding eigenmode and adoint eigenmode respectively. Then there exists no smooth and bounded solution $\phi$ to 
$$\begin{cases}\orr_{c,\nu}\phi=\psi;\\
\phi(0)=0;\\
\phi'(0)=\nu^{\gamma}\phi''(0),
\end{cases}
$$
when $\psi=\bar{\phi}_{c,\nu}$ and when $\psi=\phi^*_{c,\nu}$. The same result holds for $\orr^*_{c,\nu}$ with $\psi=\bar{\phi}^*_{c,\nu}$ and $\psi=\phi_{c,\nu}$.
\end{enumerate}
\end{theorem}
Two questions remain unclear. The first is which values the integer $\kappa$ can actually take; this might depend on the shear flow $U_s$ as well as the Rayleigh eigenvalue $c_0$, and $\alpha$. The second is whether the Orr-Sommerfeld eigenmode $\phi_{c_{\nu},\nu}$ belongs to the image $\Img \orr_{c_{\nu},\nu}$ of the Orr-Sommerfeld operator. This is not true as long as $\phi_{c_0}\notin \Im \ray_{c_0}$, where $c_{\nu}\to c_0$ as $\nu \to 0$ and $c_0$ is a Rayleigh eigenvalue. Furthermore, the latter is equivalent to
$$ \int_0^{\infty}\frac{\phi_{c_0}^2}{U_s-c_0}\neq 0.$$
In \cite{greniernguyen}, the authors make this assumption in order to justify the existence of a projector to $\Img \orr_{c_{\nu},\nu}$ by adding a multiple of a non-trivial element of $\ker \orr_{c_{\nu},\nu}$.

The adjoint operators considered in this paper prove therefore useful because they can allow us to deduces properties about the kernels and images of the original operators. The rigorous proof of the existence of Orr-Sommerfeld eigenvalues near Rayleigh eigenvalues is also new in the literature to the authors' knowledge. These properties will be employed in an upcoming work by the authors on the instability of Navier-Stokes with the Navier boundary condition, which will improve the result obtained in \cite{prevpaper}. We point out that while we wished to introduce the adjoint operators in this paper, the entire construction and techniques also extend to the original operators - see Appendix \ref{original}.

\section{The adjoint Rayleigh operator}\label{sec2}
In this section, we consider the previously introduced Rayleigh operator as a densely defined operator in $L^2(\mathbb{R}_+;\mathbb{C})$ with domain
\begin{align*}
    \mc{D}_{0}&:=\lll \phi \in H^{2}(\mathbb{R}_+;\mathbb{C}) : \phi(0)=0 \rrr,
\end{align*}
though of course $\ray_c\phi$ makes sense for any $\phi$ with sufficient regularity. In this Section, we assume by simplicity that $\alpha>0$ as the inviscid equations are invariant under a change of sign of $\alpha$. We introduce the adjoint operator $\ray_c^*$ which is obtained via integration by parts. We will frequently use the partial derivative ($\partial_y$) notation as we also apply these operators to two-variable functions when we introduce the Green functions.

\begin{prop}\label{rayadjdef} The Rayleigh operator has an adjoint on $L^2(\mathbb{R}_+;\mathbb{C})$, given by
\begin{align}
\ray_c^*&= \Delta_{\alpha}(U_s-\bar{c})-U_s''=(U_s-\bar{c})\Delta_{\alpha}+2U_s'\partial_y,\qquad \Delta_{\alpha}=\partial_{yy}-\alpha^2,\nonumber \end{align}
whose domain contains $\mc{D}_0$. Moreover, we have
\begin{align}\label{eq:rayadj}\ray_c^*\phi &= \frac{1}{U_s-\bar{c}} \overline{  \ray_c((U_s-c)\bar{\phi})}. \end{align}

\end{prop}
\begin{proof}Let $\phi,\psi\in\mc{D}_0$. Integrating by parts,
\begin{align*}
    \int_0^{\infty}\ray_c\phi \bar{\psi} &= \int_0^{\infty}(U_s-c)\Delta_{\alpha}\phi \bar{\psi}-\int_0^{\infty}U_s'' \phi \bar{\psi}=\int_0^{\infty}\phi \Delta_{\alpha}((U_s-c)\bar{\psi})-\int_0^{\infty}\phi U_s''\bar{\psi}= \\
    &=\int_0^{\infty}\phi \overline{ \l \Delta_{\alpha}(U_s-\bar{c})-U_s''\r \psi},
\end{align*}
where we could integrate by parts thanks to the boundary condition $\phi(0)=\psi(0)=0$.

Let us check \eqref{eq:rayadj}. We have
\begin{align*}
  \frac{1}{U_s-\bar{c}}\overline{  \ray_c((U_s-c)\bar{\phi})}&=\frac{1}{U_s-\bar{c}}\overline{ (U_s-c)\Delta_{\alpha}((U_s-c)\bar{\phi})-(U_s-c)U_s''\bar{\phi}}\\ 
  &= \Delta_{\alpha}((U_s-\bar{c})\phi)-U_s''\phi=\ray_c^*\phi.
\end{align*}
\end{proof}
We can also define eigenvalues and eigenmodes for the adjoint operator.
\begin{deff}\label{rayeig}
A complex number $c=c_0\in \mathbb{C}$ is called an \emph{eigenvalue} for $\ray_c$ (respectively, $\ray_c^*$) when there exists $\phi \in \mc{D}_0$ with $\ray_{c_0}\phi=0$ (resp. $\ray_{c_0}^*\phi=0$), i.e. when $\ker \ray_{c_0}\neq \lll 0\rrr$ (resp. $\ker \ray_{c_0}^*\neq \lll 0 \rrr$).
\end{deff}
However, from \eqref{eq:rayadj} we immediately deduce that \emph{$c$ is an eigenvalue for $\ray_c$ if and only if $c$ is an eigenvalue for its adjoint $\ray^*_c$}. From now on, we will simply refer to them as \emph{Rayleigh eigenvalues}, without distinguishing between the original and adjoint operator. 
\subsection{Properties of the fundamental solutions}\label{fund}

In this section, we construct two solutions of the adjoint Rayleigh equation $\phi^{*,\pm}_c$ with respective behavior at infinity $Ce^{\pm \alpha y}$, for some $C>0$. These solutions could also be obtained by setting
$$\phi^{*,\pm}_c=\frac{\bar{\phi}^{\pm}_c}{U_s-\bar{c}},$$
where $\phi^{\pm}_c$ is a non-trivial decaying or growing solution of Rayleigh. Their existence and asymptotic properties are often stated in the literature, but we were not able to find a rigorous proof. The methods of this section do not require the prior existence of these Rayleigh solutions. In fact, they can be easily adapted to prove their existence and similar properties to the one obtained for the adjoint equation (see Appendix \ref{original}).

For all $\bar{y}\geq 0$, $\eta \in \mathbb{R}$ and $k\in \mathbb{N}$, we introduce the space of exponentially bounded functions with $k$ derivatives on $[\bar{y},\infty)$:
\begin{equation}\label{eq:space}W^{\infty}_{k,\eta;\bar{y}}:=\lll \phi\in C^k(\mathbb{R}_+) : \|\phi\|_{k,\eta}:=\max_{0\leq j\leq k}\sup_{y\geq \bar{y}}e^{\eta y}|\partial^j\phi(y)|<\infty\rrr .\end{equation}
We simplify the notation by dropping the first or third parameter when it is equal to $0$, e.g.
$W^{\infty}_{k,\eta}:=W^{\infty}_{k,\eta;0}$ and $L^{\infty}_{\eta;\bar{y}}:=W^{\infty}_{0,\eta;\bar{y}}$, and analogously for the respective norms. We denote the standard $L^{\infty}$-based Sobolev space of order $n$ with $W^{n,\infty}$.

\begin{lemma}\label{exp}Let $\alpha >0$, $\eta_2>-\alpha$ and $\eta_1<\eta_2$, with $\eta_1\leq \alpha$. Then
\begin{equation}\sup_{y\geq 0} e^{\eta_1y} \int_0^{\infty}e^{-\alpha|x-y|}e^{-\eta_2x}\,\mathrm{d}x\leq \begin{cases}\frac{2\alpha }{(\alpha-\eta_2)(\alpha+\eta_2)} & \eta_2 \in (-\alpha, \alpha);\\
\frac{\alpha+\eta_1}{2\alpha(\alpha-\eta_1)}& \eta_2= \alpha;\\ 
 \frac{1}{\eta_2-\alpha} & \eta_2> \alpha \geq \eta_1.\end{cases}
\end{equation}
\end{lemma}
\begin{proof}
For all $\eta_2>-\alpha$, we have
\begin{align*}
\int_0^{\infty}e^{-\alpha|y-x|}e^{-\eta_2 x}\,\mathrm{d}x= \begin{cases}e^{-\alpha y}\frac{2\alpha e^{(\alpha-\eta_2)y}-(\eta_2+\alpha)}{(\alpha+\eta_2)(\alpha-\eta_2)} & \eta_2\neq \alpha;\\
e^{-\alpha y} \frac{1+2\alpha y}{2\alpha} & \eta_2=\alpha.
\end{cases}
\end{align*}
We can then obtain the result by neglecting all terms with a negative sign and using the fact that the supremum over $\mathbb{R}_+$ of a negative exponential is equal to $1$.
\end{proof}
Let $\alpha>0$. For all $\phi\in L^{\infty}_{-\alpha;\bar{y}}$, $\bar{y}\geq 0$ we may define

\begin{equation}\label{eq:top}
   T_{\alpha;\bar{y}}\phi(y):=\frac{1}{2\alpha(U_s(y)-\bar{c})}\int_{\bar{y}}^{\infty}e^{-\alpha |y-x|}U_s''(x)\phi(x)\,\mathrm{d}x.
\end{equation}
Note that $e^{-\alpha|y-x|}$ is a Green function for the $\Delta_{\alpha}$ operator according to Proposition \ref{generalgreen}, and therefore we deduce that
$\Delta_{\alpha}(U_s-  \bar{c})T_{\alpha;\bar{y}}\phi=U_s''\phi$ on $(\bar{y},+\infty)$. Moreover, since $U_s$ is smooth, $T_{\alpha;\bar{y}}$ maps smooth functions to smooth functions. 

The following Lemma will play a crucial role in constructing solutions of both the adjoint Rayleigh and Orr-Sommerfeld.  

\begin{lemma}\label{Top} Let $\bar{y}\geq 0$, $n\geq 0$ and let $\eta_0$ be the constant from \eqref{eq:eta}. The operator $T_{\alpha;\bar{y}}$ maps $W^{\infty}_{n,\eta;\bar{y}}$ to itself if (and only if) $\eta \in (-\alpha-\eta_0,\alpha]$. Moreover, there exists $y_n\geq 0$ such that $T_{\alpha;\bar{y}}$ is a contraction on $W^{\infty}_{n,\pm \alpha;\bar{y}}$ for $\bar{y}\geq y_n$.
\end{lemma}
\begin{proof}We can assume $\eta_0<2\alpha$. Taking the $n$-th derivative of \eqref{eq:top}, we have
\begin{equation}\label{eq:talpha}
    \|\partial^n T_{\alpha;\bar{y}}\phi\|_{\eta;\bar{y}}\leq  \max_{0\leq j\leq n}\|\partial^j \phi\|_{\eta;\bar{y}} C_{k}\alpha^{n-1}\| U_s''\|_{n,\eta;\bar{y}}\sup_{y\geq \bar{y}}e^{\eta y}\int_{y_0}^{\infty}e^{-\alpha |x-y|}e^{-(\eta_0+\eta)x}\,\mathrm{d}x,
\end{equation}
for some constant $C_{n}>0$.

By Lemma \ref{exp}, if $\eta \in (-\alpha-\eta_0,\alpha]$ then the supremum is finite and bounded by 
$$ \begin{cases}2\alpha(\alpha-\eta_0-\eta)^{-1}(\alpha+\eta_0+\eta)^{-1}&\quad \text{ if }\eta+\eta_0<\alpha\\ 
(\alpha+\eta)(2\alpha(\alpha-\eta))^{-1} & \quad\text{if }\eta+\eta_0=\alpha\\ 
(\eta+\eta_0-\alpha)^{-1} & \quad\text{ if }\eta + \eta_0>\alpha.
\end{cases}
$$
From \eqref{eq:talpha}, we deduce the mapping property. 
\begin{comment}To prove the first contraction property, notice
$$\|T_{\alpha}\|_{k,\alpha}\leq \frac{C_0\|U_s''\|_{k,\alpha;0}}{\alpha\eta_0},\qquad \|T_{\alpha}\|_{k,-\alpha}\leq \frac{2C_0 \|U_s''\|_{k,-\alpha;0}}{\eta_0(2\alpha-\eta_0)}, $$
both of which vanish as $\alpha \to \infty$. 
\end{comment}
Since $\partial^k U_s''\in L^{\infty}_{\eta_0}$ for all $k\geq 0 $, then for each $n\geq 0$ and for each given $\tilde{C}_n>0$ there exists $y_n\geq 0$ such that
\begin{equation}\label{eq:yous}
   \max_{0\leq k\leq n} C_{k}\alpha^{k-1}\cdot\|U_s''\|_{n,\eta_0/2;y_n} \leq \tilde{C}_{n},
\end{equation} By \eqref{eq:yous} and \eqref{eq:talpha}, for any $\bar{y}\geq y_n$ we have

$$\|T_{\alpha;\bar{y}}(\phi)\|_{n,-\alpha;y_n}\leq \tilde{C}_{n}\frac{2\alpha}{(2\alpha-\eta_0/2)\eta_0/2}\|\phi\|_{n,-\alpha;y_n},\qquad \|T_{\alpha;\bar{y}}(\phi)\|_{L^{\infty}_{n,\alpha;y_n}}\leq \tilde{C}_n\frac{2}{\eta_0}\|\phi\|_{k,\alpha;y_n}.$$
In both cases, we can choose $\tilde{C}_n$ small enough to achieve the desired result.
\end{proof}
Notice that the adjoint Rayleigh equation is a second order linear ODE with coefficients that are smooth and bounded in $y\geq 0$. Therefore, if we fix a pair $(\phi(\bar{y}),\phi'(\bar{y}))$ for some $\bar{y}\geq 0$, there exists a unique smooth solution defined over $[0,\infty)$ satisfying this condition. This will be required in the proof of the following result.
\begin{prop}\label{initial}For all $c\in \mathbb{C}$ there exist two smooth solutions $\phi_c^{*,-}$ and $\phi_c^{*,+}$ of \eqref{eq:raystar} such that for all $n\in \mathbb{N}$ we have
\begin{equation}\label{eq:asymptt}\partial^n \phi_c^{*,\pm}(y)\sim C_ne^{\pm \alpha y},\qquad y\to \infty.\end{equation}
\end{prop}
\begin{proof} Fix $n\in \mathbb{N}$.
Let $y_n\geq 0$ be the value given by  Lemma \ref{Top}. Define the recursion
\begin{equation*}
    \begin{cases}\psi^{0,\pm}(y):=\frac{1}{U_s-\bar{c}}e^{\pm \alpha y};\\
    \psi^{j,\pm}:=T_{\alpha;y_n}(\psi^{j-1,\pm}) & j\geq 1.
    \end{cases}
\end{equation*}
Since $T_{\alpha;y_n}$ is a contraction over $W^{\infty}_{n,\pm \alpha;y_n}$, the series of $\psi^{j,\pm}$  converges to some $\phi_c^{*,\pm} \in W^{\infty}_{n,\mp \alpha;y_n}$. Because $\Delta_{\alpha}((U_s-\bar{c})T_{\alpha;y_n}\phi)=U_s''\phi$, we have $\Delta_{\alpha}((U_s-\bar{c})\psi^{j,\pm}) \in W^{\infty}_{n,\mp \alpha;y_n}$ and the series of $\Delta_{\alpha}((U_s-\bar{c})\psi^{j,\pm})$ converges to  $\Delta_{\alpha}((U_s-\bar{c})\phi_c^{*,\pm})\in W^{\infty}_{n,\mp \alpha;y_n}$. Hence $\partial_y^{k}\phi_c^{*,\pm}\sim C_{k} e^{\pm \alpha y}$ for all $k \leq n$.
Furthermore, we have
$$\ray^*_c(\phi_c^{*,\pm})=\sum_{j=1}^{\infty}\l \Delta_{\alpha}((U_s-\bar{c})\psi^{j,\pm})-U_s''\psi^{j,\pm}\r =\Delta_{\alpha}((U_s-\bar{c})\psi^{0,\pm})=0.$$
Thus we have constructed a solution $\phi_c^{*,\pm}$ satisfying the required properties on $[y_n,+\infty)$. We can now uniquely continue this to a smooth solution of the adjoint Rayleigh equation in $[0,\infty)$. This proves \eqref{eq:asymptt} up to order $n$. Since $n$ is arbitrary and the decaying or growing solutions constructed must coincide across the different values of $n$ by uniqueness, this unique solution satisfies \eqref{eq:asymptt} for any $n\in \mathbb{N}$.
\end{proof}

Let
\begin{equation}\label{eq:inviscid}
    J_c^*(x):=\det \begin{pmatrix}\phi_c^{*,-}(x) & \phi_c^{*,+}(x) \\ 
    \partial_y \phi_c^{*,-}(x) & \partial_y \phi_c^{*,+}(x)\end{pmatrix}.
\end{equation}
Note that if $J_c^*(x)=0$ for some $x$, then this would imply that $\phi_c^{*,-}$ and $\phi_c^{*,+}$ define the same solution of the adjoint Rayleigh equation, which is impossible. Therefore $J_c^*(x)\neq 0$. Moreover, we have $\lim_{x\to+\infty}J_c^*(x)=2C\alpha(\lim_{y\to \infty}U_s(y)-\bar{c})^{-1}$, for some $C>0$ which depends on the normalization of the fundamental solutions. As a result, $J_c^*(x)$ is bounded away from $0$. We can then define the Green function
\begin{equation}\label{eq:greenray}
    G_c^I(x,y):=-\frac{1}{(U_s(x)-\bar{c})J_c^*(x)}\begin{cases}\phi_c^{*,-}(x)\phi_c^{*,+}(y) & y<x;\\ 
    \phi_c^{*,+}(x)\phi_c^{*,-}(y) & y>x.
    \end{cases}
\end{equation}
Note that $G_c^I$ has not been defined to satisfy any specific boundary condition; those will be fixed later. It is however invariant with respect to renormalizations of the fundamental solutions. It is continuous, with a piecewise continuous and integrable derivative satisfying
$$|\partial_y G_{c}^I(x,y)|\lesssim \frac{1}{|(U_s(x)-\bar{c})J_c^*(x)|}e^{-\alpha |y-x|}.$$
The above can be considered as a Green function in the sense that for $\phi$ continuous and bounded over $\mathbb{R}_+$, then \begin{equation} \label{eq:delta}\ray_c^*\int_0^{\infty} G_c^I(x,y)\phi(x)\,\mathrm{d}x=\phi(y), \qquad \forall y \geq 0.\end{equation}
To prove this, one notices that $G_c^I(x,y)\in C^0(\mathbb{R}^2_+)$ is smooth outside of the line $x=y$, and at $x=y$ we have the jump
$$\partial_2 G_c^I(y_-,y)-\partial_2 G_c^I(y_+,y)=\frac{1}{U_s(y)-\bar{c}},\qquad \forall y\geq 0.$$
In other words, $G_c^I$ satisfies the assumptions of Proposition \ref{generalgreen} from the Appendix, which guarantees \eqref{eq:delta}.

 Observe now the identity
\begin{equation}\label{eq:rayid}
    \ray^*_c = \ray^*_{c_0}-(\bar{c}-\bar{c}_0)\Delta_{\alpha}.
\end{equation}
Define the operator
\begin{equation}\label{eq:sop} S_{c;c_0}\psi(y):=-2(\bar{c}-\bar{c}_0)\int_0^{\infty}G_{c_0}^I(x,y)\frac{U_s'(x)}{U_s(x)-\bar{c}}\psi'(x)\,\mathrm{d}x,\end{equation}
which by \eqref{eq:delta} satisfies
\begin{equation}\label{eq:seq}\ray^*_{c_0} S_{c;c_0}\psi=-2(\bar{c}-\bar{c}_0)\frac{U_s'}{U_s-\bar{c}}\psi'.\end{equation}
Integrating by parts, using the fact that $U_s'(0)=0$, we see that
\begin{align*}
S_{c;c_0}\psi(y)=&2(\bar{c}-\bar{c}_0)\int_0^{\infty}\partial_x \l G_{c_0}^I(x,y)\frac{U_s'(x)}{U_s(x)-\bar{c}}\r \psi(x)\,\mathrm{d}x;\\
\partial_yS_{c;c_0}\psi(y)=&-2(\bar{c}-\bar{c}_0)\int_0^{\infty}\partial_y G_{c_0}^I(x,y)\frac{U_s'(x)}{U_s(x)-\bar{c}}\psi'(x)\,\mathrm{d}x;\\
\partial_{yy}S_{c;c_0}\psi(y)=&-2(\bar{c}-\bar{c}_0) \int_0^{\infty}\l 2\frac{U_s'(y)}{U_s(y)-\bar{c}}\partial_y G_{c_0}^I(x,y)+\alpha^2 G_{c_0}^I(x,y)\r\frac{U_s'(x)}{U_s(x)-\bar{c}}\psi'(x)\,\mathrm{d}x\\ 
&-2(\bar{c}-\bar{c}_0)\int_0^{\infty}\frac{U_s'(y)}{(U_s(y)-\bar{c})^2}\psi'(y)\,\mathrm{d}x.
\end{align*}
The above expression can in fact be differentiated arbitrarily many times by applying \eqref{eq:delta} repeatedly (also see the proof of Proposition \ref{generalgreen}). We thus obtain that for all $n\in \mathbb{N}$ there exists a constant $C_n>0$ such that
\begin{equation}\label{eq:sest}\max_{0\leq k\leq n}\|\partial_y^k S_{c;c_0}(\psi)\|_{L^{\infty}}\leq C_n\alpha^{2n} |c-c_0|\max_{k\leq n-1\lor 0}\|\partial^k \psi\|_{L^{\infty}}.\end{equation}
Here $C_n$ depends on the derivatives of $U_s$ up to order $n$, which may grow very quickly.

The next proposition establishes that, under an appropriate normalization, the decaying solution $\phi_c^{*,-}$ approaches $0$ polynomially as $c$ approaches a Rayleigh eigenvalue $c_0$. 

\begin{prop}\label{newsol}Fix $c_0\in \mathbb{C}$, $\Im c_0>0$. For all $c\neq c_0$ with $|c-c_0|$ small enough, there exists a solution $\tilde{\phi}_c^{*,-}$ of the adjoint Rayleigh equation, which is a multiple of the slow decaying solution $\phi_c^{*,-}$, with the following properties:
\begin{enumerate}
\item for $|c-c_0|$ small enough, the function $\partial_y^k\tilde{\phi}_c^{*,-}$ is analytic in $\bar{c}$ for all $k\in \mathbb{N}$;
\item  for all $n\in \mathbb{N}$ there exists $C_n> 0$ such that as $c\to c_0$ we have
\begin{equation}\label{eq:raytay}
    \|\tilde{\phi}_c^{*,-}-\phi_{c_0}^{*,-}\|_{W^{n,\infty}}\sim C_n|c-c_0|.
\end{equation}
In particular, at the origin for some integer $\kappa>0$ we have
\begin{equation}\label{eq:raytay2}\tilde{\phi}_c^{*,-}(0)\sim C_{\kappa}(\bar{c}-\bar{c}_0)^{\kappa},\end{equation}
as $c\to c_0$, where
$$C_{\kappa}:=\lim_{c\to c_0}(\bar{c}-\bar{c}_0)^{-\kappa}(S_{c;c_0})^{\kappa}(\phi_{c_0}^{*,-})(0).$$
\end{enumerate}

\end{prop}
\begin{proof}
Let $\psi^0:=\phi_{c_0}^{*,-}$, which satisfies $\ray_{c_0}^*\phi_{c_0}^{*,-}=0$. Recall the operator $S_{c;c_0}$ defined in \eqref{eq:sop} and define recursively for $j\geq 1$
$$ \psi^{j}:=S_{c;c_0}(\psi^{j-1}),$$
so that by \eqref{eq:seq} we have $\ray^*_{c_0}\psi^j = -2(\bar{c}-\bar{c}_0)\frac{U_s'}{U_s-\bar{c}}\partial_y\psi^{j-1}$.
Moreover, by \eqref{eq:sest} for all $n\in\mathbb{N}$ there exists a constant $C_n=C_n(c_0;c)$, smooth in $c_0$ and $c$ (since $|U_s-\bar{c}|$ is bounded away from zero), such that 
\begin{equation}\label{eq:estpar}\max_{0\leq k\leq n}\|\partial_y^k\psi^j\|_{L^{\infty}}\leq C_n^j \alpha^{2nj} |c-c_0|^j\max_{0 \leq k\leq n-j\lor 0} \|\partial_y^k\phi_{c_0}^{*,-}\|_{L^{\infty}}. \end{equation}
Here we have used the notation $a \lor b:= \max\lll a;b\rrr$. Similarly we will use $a\land b:= \min \lll a;b\rrr$.
Taking 
$$\tilde{\phi}_c^{*,-}:=\sum_{j=0}^{\infty}\psi^j, $$
this series then converges with its first $n$ derivatives for $|c-c_0|< (C_n\alpha^{2n})^{-1}$. In particular, $\partial_y^k\tilde{\phi}_c^{*,-}$ is analytic in $\bar{c}-\bar{c}_0$ and hence in $\bar{c}$ with radius of convergence $(C_k\alpha^{2k})^{-1}$. We then have
\begin{align*}\ray^*_c\tilde{\phi}_c^{*,-}&=\sum_{j=0}^{\infty}\ray^*_c\psi^j =\sum_{j=0}^{\infty}\l \ray^*_{c_0}\psi^j-(\bar{c}-\bar{c}_0)\Delta_{\alpha}\psi^j\r\\ 
&=(\bar{c}-\bar{c}_0)\sum_{j=0}^{\infty}\left[-2\frac{U_s'}{U_s-\bar{c}}\partial_y\psi^{j}-\Delta_{\alpha}\psi^j\right]= -\frac{\bar{c}-\bar{c}_0}{U_s-\bar{c}}\sum_{j=0}^{\infty}\ray^*_c\psi^j\\ 
&=-\frac{\bar{c}_0-\bar{c}}{U_s-\bar{c}}\ray^*_c\tilde{\phi}_c^{*,-}.\end{align*}
But for $c\neq c_0$ this equality is only possible if $\ray^*_c\tilde{\phi}_c^{*,-}=0$. The result then follows from \eqref{eq:estpar}.
\end{proof}
\begin{remark} The solutions $\tilde{\phi}_c^{*,-}$ constructed above, due to their asymptotics, must be linearly dependent with $\phi_c^{*,-}$, but the properties of the solution constructed in Proposition \ref{newsol} are independent from any renormalization. Therefore, these properties will hold for our fundamental solution $\phi_c^{*,-}$, regardless of the normalization.

\end{remark}
\begin{remark} In practice, we have not found a way to determine the value of $\kappa$ from Proposition \ref{newsol}. We know that
$$\kappa= \bar{\kappa}\iff (S_{c;c_0})^{\kappa}\phi_{c_0}^{*,-}(0)=0\;\forall \kappa=1,\dots, \bar{\kappa}.$$
For instance, from Corollary \ref{dep} it will follow that
$$\kappa >1 \iff \int_0^{\infty} \frac{{\phi_{c_0}^-(\phi_{c_0}^-)'}}{(U_s-c)^2}U_s'=0.$$
Similar to the integral from Remark \ref{eigenimage}, we have not been able to establish whether this integral is zero or non-zero for an arbitrary Rayleigh eigenvalue $c_0$.
\end{remark}
\subsection{Image and kernel}
In this Section, we wish to study the image and kernel of the original and adjoint Rayleigh operator on the domain $\mc{D}_0\subset L^2(\mathbb{R}_+)$. To this purpose, it is important to point out that all the results obtained in Section \ref{fund} also apply with minimal modifications to the original Rayleigh equation (see Section \ref{original}). In particular, one can prove similarly to Proposition \ref{initial} that the Rayleigh equation \eqref{eq:ray}  has two fundamental solution $\phi_c^{\pm}$ with $\phi_{c}^{\pm}(y)\sim e^{\pm \alpha y}$ as $y\to \infty$ (see Appendix \ref{original}). We use the symbol $\Img$ to denote the image of an operator, to differentiate it from the imaginary part of a complex number.

If $c$ is not a Rayleigh eigenvalue, then by definition the operators $\ray_c$ and $\ray^*_c$ are injective. Furthermore, as $\phi_c^-(0)\neq 0$ one may define a Green function 
$$G_c(x,y)=G_c^I(x,y)- \frac{G_c^I(x,0)}{\phi_c^-(0)}\phi_c^-(y), $$
which will satisfy the boundary condition $G_c(x,0)=0$. Therefore, any sufficiently regular function $\phi$  such that $\int_0^{\infty}G_c(x,y)\phi(x)\,\mathrm{d}x$ converges will belong to $\Im \ray^*_c$ by Proposition \ref{generalgreen}. 

The picture becomes more interesting if we let $c=c_0$ be a Rayleigh eigenvalue. The above argument does not extend to this case as $\phi_c^-(0)=0$ by definition. We start with the following elementary Lemma.
\begin{lemma}\label{image0} Let $c_0$ be a Rayleigh eigenvalue. Then
$\Delta_{\alpha}\bar{\phi}_{c_0}^-$ and $\frac{U_s''\bar{\phi}_{c_0}^-}{U_s-\bar{c}_0}$ belong to $\Img \ray_{c_0}$. Similarly, $\Delta_{\alpha}\bar{\phi}_{c_0}^{*,-}$ and $2U_s'\frac{\partial_y\bar{\phi}_{c_0}^{*,-}}{U_s-c_0}$ belong to $\Img \ray^*_{c_0}$.
\end{lemma}
\begin{proof}
    A direct computation shows that
    $$ \ray_{c_0}\bar{\phi}_{c_0}= (\bar{c}_0-c_0)\Delta_{\alpha}\bar{\phi}_{c_0}=U_s''\frac{\bar{c}_0-c_0}{U_s-\bar{c_0}}\bar{\phi}_{c_0},$$
    and $\bar{\phi}_{c_0}\in \mc{D}_0$ because $\bar{\phi}_{c_0}(0)=0$. For the adjoint operator, the proof is similar.
\end{proof}

From \eqref{eq:rayadj}, we conclude that the two fundamental solutions to $\ray_c^*=0$ are (up to a multiplicative constant)
\begin{equation}\label{eq:fundrayadj}\phi^{*,\pm}_c= \overline{\phi_c^{\pm}(U_s-c)^{-1}}=\frac{\overline{\phi^{\pm}_c}}{U_s-\bar{c}}.\end{equation}
Recall that $U_s-\bar{c}$ cannot vanish by assumption \eqref{eq:critlay}. Therefore the fundamental solutions are well defined and enjoy the same asymptotic behavior as the original fundamental Rayleigh solutions:
$$\phi_c^{*,\pm}(y)\sim e^{\pm \alpha y},\qquad y\to \infty. $$

The following Lemma is a peculiarity of the original Rayleigh equation. It does not hold for the adjoint equation due to the presence of the additional term of order $1$ which breaks the symmetry.
\begin{lemma}\label{jconst} Consider
\begin{equation}
    J_c(x):=\det \begin{pmatrix}\phi_c^-(x) & \phi_c^+(x) \\
    \partial_y \phi_c^-(x) & \partial_y \phi_c^+(x) \\
    \end{pmatrix}.
\end{equation}
Then $J_c(x)=C\alpha$ for all $x\geq 0$, where $C\in \mathbb{C}\setminus \lll 0 \rrr$ is independent from $\alpha$. In particular, if $c_0$ is a Rayleigh eigenvalue, then 
$$\partial_y \phi_{c_0}^-(0)=\frac{C\alpha}{\phi_{c_0}^+(0)}.$$
\end{lemma}
\begin{proof}
    We have
    $$J_c'(x)= \phi_c^-\partial_{yy}\phi_c^+-\phi_c^+\partial_{yy}\phi_c^+=\phi_c^-\l \alpha^2 + \frac{U_s''}{U_s-c}\r\phi_c^+-\phi_c^+\l \alpha^2+\frac{U_s''}{U_s-c}\r\phi_c^-=0.$$
    Therefore $J_c$ is constant, and it must be non-zero or it would violate the uniqueness of solutions for the initial value problem for the Rayleigh equation. In fact, we have
    $$J_c(x)=\lim_{x\to \infty}J_c(x)= C\alpha,$$
    for some $C\in \mathbb{C}\setminus \lll 0\rrr$ independent from $\alpha$. If $c=c_0$ is a Rayleigh eigenvalue, then recall that $\phi_{c_0}^-(0)=0$ and $\phi_{c_0}^+(0)\neq 0$.
\end{proof}
\begin{lemma}\label{jconst2}Let $J_c^*(x)$ be the corresponding matrix for the adjoint operator defined in \eqref{eq:inviscid}. Then we have
\begin{equation}\bar{J_c^*}(x)=\frac{1}{(U_s(x)-c)^2}J_c(x). 
\end{equation}
In particular, by Lemma \ref{jconst} we have that $(U_s(x)-c)^2\bar{J_c^*}(x)=C\alpha$, where $C$ is the same constant from the Lemma.
\end{lemma}
\begin{proof}
    Recalling that $\bar{\phi}_c^{*,\pm}=\frac{{\phi}_c^{\pm}}{U_s-{c}}$, we have 
    \begin{align*}(U_s-c)^3\bar{J_c^*}(x)=\det \begin{pmatrix} \phi_c^- & \phi_c^+ \\ 
    \partial_y \phi_c^- (U_s-c)-U_s'\phi_c^-& \partial_y \phi_c^+(U_s-c)-U_s'\phi_c^+
    \end{pmatrix}(x)=(U_s(x)-c)J_c(x).
\end{align*}
\end{proof}
We conclude the following result, which proves that the conjugate of the Green function $G_{c}^I$ fron \eqref{eq:greenray} is a multiple of the eigenmode when $c=c_0$ is a Rayleigh eigenvalue.
\begin{corollary}\label{dep}Let $\phi_{c_0}^-$ be an eigenmode for the original Rayleigh equation. Then 
$\overline{G_{c_0}^I(x,0)}$ and $\phi_{c_0}^-(x)$ are linearly dependent.
\end{corollary}
\begin{proof}
    Indeed from \eqref{eq:greenray} and \eqref{eq:fundrayadj} we have
$$\overline{G_{c_0}^I(x,0)}=\frac{\phi_{c_0}^{*,+}(0)}{(U_s(x)-{c_0})\bar{J_{c_0}^*}(x)}\frac{\phi_{c_0}^-(x)}{U_s(x)-c_0}.$$
But Lemma  \ref{jconst2} tells us that  $(U_s(x)-c_0)^2\bar{J_{c_0}^*}(x)$ is a non-zero constant in $x$, so the factor multiplying $\phi_{c_0}^-(x)$ is a non-zero constant.
\end{proof}
Another consequence of \eqref{eq:fundrayadj} is the following statement on the image of the Rayleigh operator. 
\begin{prop}\label{rayim}Assume that $c=c_0$ is an eigenvalue for Rayleigh, with $\Im c_0\neq 0$ and eigenfunction $\phi_{c_0}^-$. Then $ \Img \ray_{c_0}\subset (\ker \ray^*_{c_0})^{\perp}$. Furthermore, suppose that $f:\mathbb{R}_+\to \mathbb{C}$
is such that either $\Re(f(U_s-\bar{{c_0}}))$ or $\Im (f(U_s-\bar{{c_0}}))$ has a fixed sign and is not identically zero, and $f{\phi}_{c_0}^-\in L^2(\mathbb{R}_+)$.  Then $f \bar{\phi}_{c_0}^-\not\in \Img \ray_{c_0}$. 

Similarly, for the adjoint we have $\Img \ray^*_{c_0} \subset \l \ker \ray_{c_0}\r^{\perp}$ and $g\bar{\phi}_{c_0}^{*,-}\notin \Img \ray^*_{c_0}$ for all $g:\mathbb{R}_+\to \mathbb{C}$ such that either $\Re (g(U_s-{c_0}))$ or $\Im (g(U_s-{c_0}))$ has a fixed sign, $ g\not \equiv 0$, and $g\bar{\phi}_{c_0}^{*,-}\in L^2(\mathbb{R}_+)$.
\end{prop}
\begin{proof}We prove the result for $\ray_{c_0}$, and for $\ray^*_{c_0}$ the proof is identical. From \eqref{eq:rayadj}, we infer that
$$\ray_{c_0}^*\phi = 0 \iff \ray_{c_0}((U_s-{c_0})\bar{\phi})=0\iff \phi \in \lpl \frac{\bar{\phi}_{{c_0}}^-}{U_s-\bar{{c_0}}}\rer. $$
Hence, $\lpl \frac{\bar{\phi}_{{c_0}}^{-}}{U_s-\bar{{c_0}}}\rer =\ker \ray_{{c_0}}^*$ is a one dimensional complex vector space. We know that $\ker \ray^*_{c_0}\cap \mc{D}_0\subset \Img \ray_{{c_0}}^{\perp}$. However, $\ker \ray^*_{c_0}$ is either trivial or, if $c$ is an eigenvalue, spanned by  $\overline{\phi_{c_0}^{-}(U_s-{c_0})^{-1}}$, which belongs to $\mc{D}_0$. Hence $\ker \ray^*_{c_0} \subset \Im\ray_{c_0}^{\perp} $ and therefore 
$$\Img \ray_{c_0}\subset \l \ker \ray^*_{c_0}\r^{\perp}=\lpl \overline{\phi_{c_0}^{-}(U_s-{c_0})^{-1}} \rer^{\perp},$$
implying that the image has codimension  at least one. Now let $f$ satisfying the above assumptions, and let us prove that $f\bar{\phi}_{c_0}^{-}\notin \Img \ray_{c_0}$.  We have
\begin{align*} \Im \lpl f\bar{\phi}_{c_0}^-,\frac{\bar{\phi}_{c_0}^-}{U_s-\bar{{c_0}}}\rer&= \Im \int_0^{\infty}\frac{|\phi_{c_0}^{-}|^2f}{U_s-{c_0}}=\int_0^{\infty}\frac{|\phi_{c_0}^{-}|^2\Im (f(U_s-\bar{{c_0}}))}{|U_s-{c_0}|^2};\\ \Re \lpl f\bar{\phi}_{c_0}^-,\frac{\bar{\phi}_{c_0}^-}{U_s-\bar{{c_0}}}\rer&= \Re \int_0^{\infty}\frac{|\phi_{c_0}^{-}|^2f}{U_s-{c_0}}= \int_0^{\infty}\frac{|\phi_{c_0}^{-}|^2\Re (f(U_s-\bar{{c_0}}))}{|U_s-{c_0}|^2}.\end{align*}
Thus if either the real or imaginary part of $f(U_s-\bar{{c_0}})$ has a fixed sign, then $\lpl f\bar{\phi}_{c_0}^-,\frac{\bar{\phi}_{c_0}^-}{U_s-\bar{{c_0}}}\rer\neq 0 $. In either case, $f\bar{\phi}_{c_0}^{-}\notin (\ker \ray^*_{c_0})^{\perp}\subset \overline{\Img \ray_{c_0}}$.

\end{proof}
\begin{remark} Some examples of functions $f$ satisfying the assumptions of Proposition \ref{rayim} include $f=1$ and $f=h( U_s-{c_0})$ for any real valued function $h$ with a fixed sign (and not identically zero). One may also choose $f=\frac{h}{U_s-\bar{{c_0}}}$, which for $h\equiv 1$ provides an alternative proof that $\phi_{c_0}^{*,-}\notin \Img \ray_{c_0}$. Note that we cannot choose $h=U_s''$, as this function will not have a fixed sign if $\Im c \neq 0$, as a consequence of Rayleigh's criterion for the instability of the linearized Euler equations. Furthermore, this would contradict Lemma \ref{image0}.
\end{remark}

In Proposition \ref{rayim} we were only able to prove  the image of the adjoint Rayleigh operator is contained to the orthogonal of the kernel of the Rayleigh operator.  The Green function we defined in \eqref{eq:greenray} can be used to prove the equality of the two subspaces. This also gives us a quantitative criterion to check whether any function belongs to the image of the (adjoint) Rayleigh operator.
\begin{prop}\label{imageray}
 Let $c_0$ be a Rayleigh eigenvalue and let $\psi \in C^{1}\cap L^2(\mathbb{R}_+)$. Then $\psi \in \Img \ray^*_{c_0}$ if and only if
\begin{equation}\label{eq:orthog}\int_0^{\infty}G_{c_0}^I(x,0)\psi(x)\,\mathrm{d}x=0,
\end{equation}   
which by Corollary \ref{dep} is equivalent to
\begin{equation}\label{eq:orthog2}
\int_0^{\infty}\bar{\phi}_{c_0}(x)\psi(x)\,\mathrm{d}x=0.
\end{equation}
In particular, $\Img \ray^*_{c_0}$ has codimension one in $L^2(\mathbb{R}_+)$, and  $\Img \ray^*_{c_0}=(\ker \ray_{c_0})^{\perp}.$
\end{prop}
\begin{proof} Suppose that \eqref{eq:orthog} is satisfied. Then letting $\phi(y):=\int_0^{\infty}G_{c_0}^I(x,y)\psi(x)\,\mathrm{d}x$, we have that $\phi\in \mc{D}_{0}$ and $\ray^*_{c_0}\phi=\psi$ by Proposition \ref{generalgreen}. Therefore, $\psi\in \Img \ray^*_{c_0}$. Thus $\Img \ray^*_{c_0}$ has codimension at most one (since $C^{1}\cap L^2$ is dense in $L^2$). But we already know by Proposition \ref{rayim} that the codimension is at least one, so the proof is complete.

Conversely, suppose $\psi=\ray^*_{c_0}\phi$, for some $\phi \in \mc{D}_0$. Then
$$\int_0^{\infty}G_{c_0}^I(x,0)\ray^*_{c_0}\phi(x)\,\mathrm{d}x=-\int_0^{\infty}\overline{\ray_{c_0}\overline{G_{c_0}^I(x,0)}}\phi(x)\,\mathrm{d}x,$$
where we could integrate by parts because $\phi(0)=0$ by assumption and $\lim_{x\to 0^+}G_{c_0}^I(x,0)=0$. But Corollary \ref{dep} tells us that $\overline{G_{c_0}^I(x,0)}\in \ker \ray_{c_0}$, so the above integral is zero. 
\end{proof}

\begin{remark}\label{eigenimage} An interesting question is whether or not the (adjoint) eigenmode belongs to the image of the (adjoint) Rayleigh operator. For instance, according to Proposition \ref{imageray} we have  $\phi_{c_0}^-\in \Img \ray_{c_0}$ if and only if
$$ \int_0^{\infty}\frac{(\phi_{c_0}^-)^2}{U_s-c_0}=0,$$
but we were unable to find evidence that this is the case for any eigenvalue $c_0$.
\end{remark}

\subsection{Estimates for the boundary value problem}
The main result of this Section, Lemma \ref{raybounds}, only holds when $c$ is not a Rayleigh eigenvalue. For this reason, in this Section we consider $c\in \mathbb{C}$ near (but different to) a Rayleigh eigenvalues $c_0$ to see how the estimates are affected by this proximity. This result can be employed to construct slow solutions to the adjoint Orr-Sommerfeld equation outside of Rayleigh eigenvalues, but we will later introduce more powerful techniques that can bypass this issue.

We start by constructing the Green function for the adjoint Rayleigh operator with boundary condition $G_c(x,0)=0$, defined as
\begin{equation}\label{eq:greenrayb}
    G_c(x,y)=G_c^I(x,y)-\frac{G_c^I(x,0)}{\phi_c^{*,-}(0)}\phi_c^{*,-}(y)= G_c^I(x,y)-\frac{\phi_c^{*,-}(x)\phi_c^{*,+}(0)}{(U_s(x)-\bar{c})J_c^*(x)\phi^{*,-}_c(0)}\phi_c^{*,-}(y).
\end{equation}
This Green function is therefore derived from $G_c^I$ by adding a smooth and decaying solution of the adjoint Rayleigh equation. This implies that $G_c$ still satisfies the assumptions Proposition \ref{generalgreen}. However, if $c$ is close to a Rayleigh eigenvalue $c_0$ then the factor of $(\phi_c^{*,-}(0))^{-1}\sim (c-c_0)^{-1}$ implies that $G_c$ blows up as $c\to c_0$. We also have 
\begin{equation}\label{eq:firstderr}|G_c(x,y)|\leq \frac{C_G}{|\alpha(c-c_0)|}e^{-\alpha|y-x|},\qquad |\partial_yG_c(x,y)|\leq \frac{C_G}{|c-c_0|} e^{-\alpha|y-x|}.\end{equation}
\begin{prop}\label{raybounds}Let $\eta_2>-\alpha, \eta_2\neq \alpha$. Let $\psi\in W^{\infty}_{n,\eta_2}$. Then there exists a constant $C_n$ such that for all $|c-c_0|$ small enough and for all $\eta_1\in [\eta_2-\eta_0,\eta_2)$, $\eta_1\leq \alpha$ we have that the smooth solution $\phi$ to the equation 
\begin{equation}\label{eq:nonhom}\begin{cases}\ray^*_c(\phi)=\psi \\
\phi(0)=\phi(\infty)=0;
\end{cases}
\end{equation}
satisfies the estimates
\begin{align}\label{eq:recurs0}
     \|\phi\|_{n,\eta_1}\leq \frac{C_{n}\alpha^{2n}}{|\eta_2^2-\alpha^2|   |c-c_0|} \| \psi \|_{n-1\lor 0,\eta_2};\\\label{eq:recurs1}
     \|\Delta_\alpha \phi\|_{n,\eta_2}\leq \frac{C_n\alpha^{2n}}{|\eta_2^2-\alpha^2|   |c-c_0|} \| \psi \|_{n,\eta_2}.
\end{align}
Note that $C_{n}$ may grow quickly with $n$. 
\begin{proof}
By Proposition \ref{generalgreen}, the solution to \eqref{eq:recurs0} is given by
 $$ \phi(y)=\int_0^{\infty}G_c(x,y)\psi(x)\,\mathrm{d}x,\qquad \phi'(y)=\int_0^{\infty}\partial_y G_c(x,y)\psi(x)\,\mathrm{d}x.$$
 Since $G_c$ has a piecewise continuous derivative, for $k=0$ and $k=1$ by \eqref{eq:firstderr} we have
$$\|\phi\|_{k,\eta_1}\leq \|\psi\|_{\eta_2}\frac{C_G}{\alpha^{1-k}|\bar{c}-\bar{c}_0|}\sup_{y\geq 0}e^{\eta_1 y}\int_0^{\infty}e^{-\alpha |y-x|}e^{-\eta_2 x}\,\mathrm{d}x,$$
so the estimate \eqref{eq:recurs0} follows from Lemma \ref{exp}. 
We now proceed by induction on $k\geq 2$. From Proposition \ref{generalgreen} we know $\ray^*_c G_c=(U_s-\bar{c})\Delta_{\alpha}G_c+2U_s'G_c'=\delta_{x=y}$, so we have
\begin{align*}\partial_y^{k+1}\phi(y)&=\partial_y^{k-1}\int_0^{\infty}\partial_{yy}G_c(x,y)\psi(x)\,\mathrm{d}x\\
&=\partial_y^{k-1}\left[\int_0^{\infty}\l \frac{-2U_s'(y)}{U_s(y)-\bar{c}}\partial_yG_c(x,y)+\alpha^2G_c(x,y)\r\psi(x)\,\mathrm{d}x\right]+\partial_y^{k-1}\frac{\psi(y)}{U_s(y)-\bar{c}},\\
&=\partial_y^{k-1}\l  \frac{-2U_s'(y)}{U_s(y)-\bar{c}}\phi'(y)+\alpha^2\phi(y) \r+\partial_y^{k-1}\frac{\psi(y)}{U_s(y)-\bar{c}}.
\end{align*}
Therefore, for some constant $C_{k}$ depending on the first $k$ derivatives of $U_s$ we have
$$\|\partial_y^{k+1}\phi\|_{\eta_1}\leq C_{k}\alpha^2\l  \|\phi\|_{k,\eta_1}+\|\psi\|_{k-1,\eta_1}\|\r\leq \frac{C_k\alpha^{2k}}{|\eta_2^2-\alpha^2||c-c_0|}\|\psi\|_{k,\eta_2}. $$

Taking the maximum over $0\leq k\leq n-1$ the estimate \eqref{eq:recurs0} then follows. As for \eqref{eq:recurs1}, we know that 
$(U_s-\bar{c})\Delta_{\alpha}\phi=\psi+2U_s' \phi',$ so by \eqref{eq:eta} and the previous estimate,
\begin{align*}
    \|\partial^k\Delta_{\alpha}\phi\|_{\eta_2}\leq C\l  C_k\|\partial^{k+1}\phi\|_{\eta_2-\eta_{0}}+\|\partial^k\psi\|_{\eta_2}\r\leq \frac{C_k\alpha^{2k}}{|\eta_2^2-\alpha^2||c-c_0|}\|\partial^k \psi\|_{\eta_2}.
\end{align*}
Taking the maximum over $0\leq k\leq n$ we reach \eqref{eq:recurs1}.
\end{proof}
\end{prop}
\subsection{Estimates for the initial value problem}
Let $a,b\in \mathbb{R}$, and let $\psi:\mathbb{R}_+\to \mathbb{R}$ be continuous. Consider the initial value problem for the adjoint Rayleigh equation
\begin{equation}\label{eq:ivp}
    \begin{cases}
        \ray^*_c \phi(y) =\psi(y) & y\geq 0,\\ 
        \phi(0)=a,\\ 
        \phi'(0)=b.
    \end{cases}
\end{equation}
Unlike the boundary value problem, we expect this problem to be well-posed, since a global solution to \eqref{eq:ivp} exists and is unique for every continuous $\psi$. 

As the fundamental solutions  $\phi_c^{*,\pm}$ constructed in Proposition \ref{initial} are independent, let
$$M_{c}(x):= \begin{pmatrix}\phi_c^{*,-}(x) & \phi_c^{*,+}(x) \\ 
    \partial_y \phi_c^{*,-}(x) & \partial_y \phi_c^{*,+}(x)\end{pmatrix}.$$
This matrix is invertible for all $c$ and its determinant $J_c^*(x)$ is bounded from below, as argued in Section \ref{fund}. We then have
$$\tilde{\phi} (y):=\begin{pmatrix} \phi_c^{*,-} & \phi_c^{*,+}\end{pmatrix}M_c^{-1}(0) \begin{pmatrix} a \\ b \end{pmatrix}. $$
Then $\tilde{\phi}$ satisfies \eqref{eq:ivp} with $\psi=0$. Taking the difference, $w= \phi-\tilde{\phi}$ satisfies \eqref{eq:ivp} with $w(0)=w'(0)=0$. Define the Green function $G_c^{\textup{IVP}}$ for the initial value problem by
\begin{equation}\label{eq:greenivp}G_c^{\textup{IVP}}(x,y)=   -\frac{1}{(U_s(x)-\bar{c})J_c^*(x)}\begin{cases} \phi_c^{*,+}(x)\phi_c^{*,-}(y)-\phi_c^{*,-}(x)\phi_c^{*,+}(y) & y>x \\ 
0 & y<x.
\end{cases}
\end{equation}
Notice that $G_c^{\textup{IVP}}(x,0)=\partial_y G_c^{\textup{IVP}}(x,0)=0$.
This Green function, unlike others constructed in this paper, is exponentially growing in $y$. However, it is zero for all $x>y$, thus making integration possible.

In conclusion, by Proposition \ref{generalgreen} the solution to \eqref{eq:ivp} is given by
\begin{equation}\label{eq:solivp}\phi(y)=\int_0^{\infty}G_c^{\textup{IVP}}(x,y)\psi(x)\,\mathrm{d}x+\tilde{\phi}(y).
\end{equation}
Note that the above integral converges for any continuous $\psi$ as $G_c^{\textup{IVP}}(x,y)=0$ for $x>y$.
\begin{prop}\label{ivpest}
    Let $k\in \mathbb{N}$, $\psi \in C^{k}([0,y_0])$, and let $\phi$ be the solution to \eqref{eq:ivp}, which is given by \eqref{eq:greenivp}. We then have
    \begin{equation}\label{eq:ivpest}\|\partial_y^k \phi\|_{L^{\infty}}\leq C_{k}\l \|\psi\|_{W^{k-2\lor 0}}+|(a,b)|\r,
\end{equation}
where $C_k$ is a positive number depending smoothly on $\alpha$ and $c$.
\end{prop}
\begin{proof} Recall that $|J_c^*(x)|$ and $U_s(x)-\bar{c}$ are bounded from below. The proof then follows from \eqref{eq:solivp} and Proposition \ref{generalgreen}.

\begin{comment} For $k=0$ and $k=1$, the inequality follows from \eqref{eq:solivp} since $G_c^{\textup{IVP}}(x,y)$ is continuous and has an integrable derivative $\partial_y G_c^{\textup{IVP}}(x,y)$. 

For $k\geq 2$, differentiating $k$ times by \eqref{eq:solivp} we obtain 
    $$\partial_y^k \phi(y)= \partial_y^{k-2}\int_0^{\bar{y}}\l \frac{-2U_s'(y)}{U_s(y)-\bar{c}}\partial_y G_c^{\textup{IVP}}(x,y)+\alpha^2G_c^{\textup{IVP}}(x,y)\r\psi(x)\,\mathrm{d}x+ \partial_y^{k-2}\psi(y)+\partial_y^k \tilde{\phi}(y).$$
We can further apply the rule $\partial_{yy}G_c^{\textup{IVP}}(x,y)= \frac{-2U_s'(y)}{U_s(y)-\bar{c}}+\alpha^2 G_c^{\textup{IVP}}(x,y)+\delta_y(x)$ repeatedly inside the integral to obtain \eqref{eq:ivpest} for all $k\geq 2$.
\end{comment}
\end{proof}
\begin{remark}\label{ivpext}
    Proposition \ref{ivpest} also holds, with the same proof, by applying the change of variables $y\mapsto \bar{y}-y$ to the problem \eqref{eq:ivp}, for some $\bar{y}>0$, leading to the problem
    \begin{equation}
    \begin{cases} (U_s-\bar{c})\Delta_{\alpha}\phi+2U_s'\phi'=\psi & \text{on } [0,\bar{y}]\\
    \phi( \bar{y})=a\\
    \phi'( \bar{y})=b.
\end{cases}
\end{equation}
In other words, the estimate \eqref{eq:ivpest} also holds when solving the equation backwards in the $y$ variable.
\end{remark}

\section{Approximate viscous solutions}\label{sec3}
In this section, we fix $\nu>0,\alpha \in\mathbb{Z}, c\in \mathbb{C}$ with $\Im c>0$. The Orr-Sommerfeld equation was defined in \eqref{eq:orrs}, paired with the boundary conditions \eqref{eq:orrbc}. We shall consider this operator on the domain
\begin{equation}\label{eq:domain2}\mc{D}_{\nu}:=\lll \phi\in H^4(\mathbb{R}_+): \phi(0)=0\textup{ and } \phi'(0)=\nu^{\gamma}\phi''(0)\rrr, \end{equation}though naturally $\orr_{c,\nu}\phi$ makes sense for any $\phi$ with sufficient regularity.
\begin{prop}The Orr-Sommerfeld operator has an adjoint on $L^2(\mathbb{R};\mathbb{C})$, whose domain contains $\mc{D}_{\nu}$, given by 
\begin{equation}
\label{eq:orradj}
\orr_{c,\nu}^*= (\partial_{yy}-\alpha^2)(U_s-\bar{c})-U_s''+\frac{\nu}{i\alpha} (\partial_{yy}-\alpha^2)^2.
\end{equation}
\end{prop}
\begin{proof}
The main point in obtaining the expression above just requires that we can integrate by parts
$$\int_0^{\infty}\phi {\psi}^{(4)}=\int_0^{\infty}\psi{\phi}^{(4)}. $$
We have
$$\int_0^{\infty}\phi {\psi}^{(4)} = \int_0^{\infty}\phi''{\psi}'' -\left[\phi'{\psi}^{''}\right]_0^{\infty}=\int_0^{\infty}\psi {\phi}^{(4)}-\left[\phi'{\psi}''\right]_0^{\infty}+\left[\phi''{\psi}'\right]_0^{\infty}.  $$
But since both $\phi$ and $\psi$ satisfy the same boundary condition, 
$$\phi'(0){\psi}''(0)=\nu^{\gamma}\phi'(0){\psi}'(0)=\phi''(0){\psi}'(0). $$
\end{proof}
Notice that we also have $\orr_{c,\nu}^*= \ray_{c}^*+\frac{\nu}{i\alpha} \Delta_{\alpha}^2 $, with the same boundary conditions as the original operator. 

Whilst the adjoint operators $\ray^*_c$ and $\orr^*_{c,\nu}$ were rigorously constructed on domains $\mc{D}_0$ and $\mc{D}_{\nu}$ respectively, both of these operators can be applied to any sufficiently regular function. The goal of this section is to identify two pairs of approximate solutions to \eqref{eq:orradj} which do not satisfy any particular boundary condition:
\begin{enumerate}
    \item two slow solutions $\phi^{*,s,\pm}$, constructed from the Rayleigh solutions and which follow the same asymptotics as $y\to +\infty$ ;
    \item two fast solutions $\phi^{*,f,\pm}$ obtained by neglecting $U_s''$, which is small when $y$ is large by \eqref{eq:critlay}.
    \end{enumerate}
    In line with \cite{grenieros}, we assume that $\alpha$ is not too large with respect to the viscosity:
$$ |\alpha| \ll \nu^{-1/2}.$$
This will come into play when constructing the fast solutions $\phi^{*,f,\pm}$ in Section \ref{fast}.
  The fast modes arise from considering the equation
$$ (U_s-\bar{c})\phi =-\frac{\nu}{i\alpha}\Delta_{\alpha}\phi \implies \phi^{*,f,\pm}(y)\sim e^{\pm \int \mu_* \,\mathrm{d}y},$$
  implying 
\begin{equation}\label{eq:muf} \mu_*=\mu_*(\alpha, c,\nu,y)= |\alpha|^{1/2}\nu^{-1/2}\sqrt{\alpha\nu-i(U_s(y)-\bar{c})},\end{equation}
where as a convention we take the root with positive real part.  This implies that $\phi^{*,f,\pm}$ is always rapidly increasing or decreasing respectively as $y\to\infty$. We have $\Re(\alpha \nu-i(U_s-\bar{c}))=\alpha\nu+\Im c>0$ for $\nu$ small enough, and $\Im (\alpha \nu-i(U_s-\bar{c}))=-\Re (U_s-\bar{c}).$  Note that the imaginary part might vanish at some points, but since the real part is always strictly positive, both square roots always have non-zero real parts. In polar form, we may write $\mu_*^2= \rho e^{i\theta}$ with $\theta \in (-\pi/2,\pi/2)$ and $\rho>0$, and the square root with positive real part is then always $\mu_*=\sqrt{\rho}e^{i\theta/2}$, which is an analytic of $\mu_*^2$. We conclude that $\mu_*$ is analytic with respect to $y$ and $c$.\\

Recall that for the original Orr-Sommerfeld equation we also have two fast solutions $\phi^{f,\pm}(y)\sim  e^{\pm \int \mu_*\,\mathrm{d}y}$, with $\mu=\alpha^{1/2}\nu^{-1/2}\sqrt{\alpha \nu +i(U_s-c)}$, again choosing the root with positive real part. 
We deduce the following fundamental result.

\begin{prop} The fast coefficient $\mu_*$ of the adjoint Orr-Sommerfeld equation defined in \eqref{eq:muf} is the complex conjugate of its counterpart $\mu$ for the original equation.
\end{prop}
\begin{proof}
 Since $U_s$ and $\alpha\nu$ are real, then $\alpha\nu-i(U_s-\bar{c})$ is the complex conjugate of $\alpha \nu+i(U_s-c)$. Let $\alpha\nu-i(U_s-\bar{c})=\rho e^{i\theta}$, $\theta \in (-\pi/2,\pi/2)$. Then $\mu_*= \sqrt{\rho}e^{i\theta/2}$. As for $\mu$, we have $\alpha\nu+i(U_s-c)=\rho e^{-i\theta}$, where we 
still have $-\theta\in (-\pi/2,\pi/2)$. Therefore, $\mu=\sqrt{\rho} e^{-i\theta/2}$ which is the complex conjugate of $\mu_*$.
\end{proof}
In this section, we wish to establish the following result.
\begin{prop}\label{appsol}For each $N\geq 1$ there exist four independent approximate solutions of the adjoint Orr-Sommerfeld equation $\phi^{*,s,\pm}_N$ and $\phi^{*,f,\pm}_N$ such that
\begin{equation}
   \partial_y^j \phi^{*,s,\pm}_N(y)\sim {\l 1+O(\nu)\r }{\partial_y^j} \phi^{*,\pm}_c(y),\qquad \partial_y^j\phi^{*,f,\pm}_N(y)\sim (1+O(\sqrt{\nu}))\partial_y^je^{-\int\mu_*\,\mathrm{d}y},\qquad \text{ as }\nu \to 0,\;j\in \mathbb{N},
\end{equation}
and for some $D_N>0$,
$$|\orr^*_{c,\nu}\phi^{*,\ell,\pm}_N(y)| \leq D_N\nu^{N+1}|\phi^{*,\ell,\pm}_N(y)|,\qquad y\in \mathbb{R}_+,\ell=s,f.$$
\end{prop}
\subsection{Approximate slow modes}
Similarly to the original Orr-Sommerfeld operator (\cite{grenieros}), the slow modes $\phi^{*,s,\pm}$ are constructed starting from the corresponding Rayleigh solutions $\phi^{*,\pm}$. Indeed, they follow the same asymptotics at infinity. 
Compared to Grenier's results for the Orr-Sommerfeld equation, we present an improved approach that successfully produces approximate adjoint Orr-Sommerfeld solutions for any value of $\alpha$, even when $c$ is near a Rayleigh eigenmode. These approximate solutions do not converge to an exact Orr-Sommerfeld solution, but they will be used later on to construct one.

Recall that the adjoint Rayleigh operator $\ray^*_{c}$ has two  smooth fundamental solutions $\phi_c^{*,\pm}$ given by $\overline{{\phi_{c}^-}}\l {U_s-\bar{c}}\r^{-1}$, where  $\phi_{c}^{\pm}$ is the corresponding fundamental solution of Rayleigh. We shall prove the following.
\begin{prop}\label{slowsol}There exist two smooth approximate solutions $\phi^{*,s,\pm}_N\in W^{\infty}_{n,\mp\alpha}$ for all $n\in \mathbb{N}$ to the adjoint Orr-Sommerfeld equation \eqref{eq:orradj}, in the form 
\begin{equation}\label{eq:adjslow}\partial_y^j\phi^{*,s,\pm}_N= \partial_y^j \phi^{*,\pm}_c+\sum_{k=1}^{N} \partial_y^j{\psi}_{c}^{k,\pm}\nu^k,\qquad \text{ as }\nu \to 0,j\in \mathbb{N}, \end{equation}
where $\|\partial_y^j\psi_c^{k,\pm}\|_{\infty,\mp\alpha}\leq C_j$, for some $C_j>0$ independent from $k$ and $\nu>0$.
\end{prop}
\begin{proof}
 By Lemma \ref{Top}, for $y_n$ large enough (depending on $U_s$ only), $T_{\alpha;y_n}$ is a contraction on $W^{\infty}_{n,\pm\alpha;y_n}$ for some $y_n\geq 0$ depending on $n\in \mathbb{N}$. Let $\phi_{c}^{*,\pm}$ be the respectively decaying or growing fundamental solution of the adjoint Rayleigh equation. Consider the iterative scheme
\begin{equation}\label{eq:iter0}\ray^*_{c}\l \psi_{c} ^{k,\pm}\r=-\frac{\nu}{i\alpha}\Delta_{\alpha}^2\psi_c^{k-1,\pm} \iff \psi_c^{k,\pm} - T_{\alpha;y_n}[\psi_c^{k,\pm}]=-\frac{\nu}{i\alpha}\l \frac{\Delta_{\alpha}\psi_c^{k-1,\pm}}{U_s-\bar{c}}\r.\end{equation}
Fix $n\geq 2$. Specifically, we choose the solution which on $[y_n,+\infty)$ is given by 
$$\psi_c^{k,\pm} :=\begin{cases} \phi_{c}^{*,\pm} & k=0 \\ 
- \frac{\nu}{i\alpha}\sum_{j=0}^{\infty}(T_{\alpha;y_n})^j \left[ \frac{\Delta_{\alpha}\psi_c^{k-1,\pm}}{U_s-\bar{c}}\right] & k\geq 1 \end{cases}.$$
Note that the above function does not depend on $n$ except through the interval of convergence.

We want to prove that $\psi_c^{k,\pm}$ can be extended to a smooth function on $\mathbb{R}_+$ for all $k\in \mathbb{N}$. The proof is by induction over $k\in \mathbb{N}$. The base case $k=0$ holds by Proposition \ref{initial}. For $k\geq 1$, if $\psi_c^{k-1,\pm}$ is smooth on $\mathbb{R}_+$, then the above series converges in $W^{\infty}_{n,\mp\alpha;y_n}$ since $T_{\alpha;y_n}$ is a contraction on this space by Lemma \ref{Top}. We therefore know that  $\psi_c^{k,\pm}\in C^2([y_n,+\infty))$ and solves \eqref{eq:iter0}. We now extend $\psi_c^{k,\pm}$ to $[0,\infty)$ by taking the solution $\phi$ to the backwards initial value problem 
\begin{equation}\label{eq:newivp} \begin{cases}
    \ray^*_c\phi= -\frac{\nu}{i\alpha}\Delta_{\alpha}^2\psi_c^{k-1,\pm} & \text{ on }[0,y_n),\\
    \phi(y_n)=\psi_c^{k,\pm}(y_n),\\
    \phi'(y_n)=(\psi_c^{k,\pm})'(y_n).
    \end{cases}
\end{equation}
By Proposition \ref{ivpest} (see Remark \ref{ivpext}), since $\psi_c^{k-1,\pm}$ is smooth on $\mathbb{R}_+$ then the solution $\phi$ is smooth on $[0,y_n)$. Extending $\psi_c^{k-1,\pm}$ as described above, the resulting solution, which we call $\phi_c^{k-1,\pm}$, is of class $C^1$ over $[0,\infty)$ and piecewise smooth. However, $\psi_c^{k,\pm}$ has a unique smooth extension $\tilde{\phi}_c^{k,\pm}$ to $[0,\infty)$. Then $\phi=\tilde{\phi}_c^{k,\pm}$ will also satisfy \eqref{eq:newivp}, and therefore by the uniqueness of solutions to \eqref{eq:newivp} we must have $\tilde{\phi}_c^{k,\pm}=\phi_c^{k,\pm}$. In conclusion, $\phi_c^{k,\pm}$ is smooth over $[0,\infty)$. Finally, since $\psi_c^{k-1,\pm}\in W^{\infty}_{n,\mp \alpha;y_n}$ then $\phi_c^{k-1,\pm}\in W^{\infty}_{n,\mp\alpha}$. This completes the inductive step. Note that the solutions constructed above cannot depend on $n$, as for each $n\geq 2$ there is a unique $C^n$ solution to \eqref{eq:iter0} on $[0,\infty)$.

Now, taking
$$\phi_N^{*,s,\pm}:=\sum_{k=0}^{N}\phi_c^{k,\pm}, $$
we get that $\phi^{*,s,\pm}_{N}$ is smooth on $\mathbb{R}_+$ as it is the sum of smooth functions, and by \eqref{eq:iter0} and \eqref{eq:newivp}, we have
\begin{equation*}\orr^*_{c,\nu}\phi_N^{*,s,\pm}=\sum_{k=0}^{N}\l \ray^*_{c}\phi_c^{k,\pm}+\frac{\nu}{i\alpha}\Delta_{\alpha}^2\phi_c^{k,\pm}\r=\frac{\nu}{i\alpha}\Delta_{\alpha}^2\phi_c^{N,\pm}\qquad \text{ on }\mathbb{R}_+. \end{equation*}
By \eqref{eq:iter0} on $[y_n,+\infty)$ and for all $k\in \mathbb{N}$ we have
\begin{equation}\label{eq:estlargey}\|\psi_c^{k,\pm}\|_{n,\mp\alpha;y_n}\leq \frac{\nu}{\alpha} \frac{1}{1-\|T_{\alpha;y_n}\|}\left\|\frac{\Delta_{\alpha}\psi_c^{k-1,\pm}}{U_s-\bar{c}}\right\|_{n,\mp\alpha;y_n}\leq C\frac{\nu}{\alpha} \|\Delta_{\alpha}\psi_c^{k-1,\pm}\|_{n,\mp\alpha;y_n}\leq C\frac{{\nu}^k}{\alpha^k}\|\Delta_{\alpha}^k\phi_c^{0,\pm}\|_{n,\mp\alpha;y_n}.\end{equation}
 Furthermore, Proposition \ref{ivpest} applied on the interval $[0,y_n]$ (see Remark \ref{ivpext}) yields the same estimate (with a different constant) for $W^{n,\infty}$ norms on this interval:
$$ \|\phi_c^{k,\pm}\|_{W^{n,\infty}}\leq C_n\frac{\nu}{\alpha}\|\Delta_{\alpha}^2\phi_c^{k-1,\pm}\|_{W^{n,\infty}}\leq C'_n \frac{\nu^k}{\alpha^k}\|\Delta_{\alpha}^{2k}\phi_c^{0,\pm}\|_{W^{n,\infty}}=C'_n \frac{\nu^k}{\alpha^k}\|\Delta_{\alpha}^{2k}\phi_c^{*,\pm}\|_{W^{n,\infty}}.$$
As $[0,y_n]$ is bounded, we can combine the two estimates together in $W^{\infty}_{n,\mp\alpha}$.
Hence,
$$\|\orr_{c,\nu}^*\phi_N^{*,s,\pm} \|_{\mp \alpha}\leq \frac{\nu^{N+1}}{\alpha^N}\|\Delta_{\alpha}^{N}\phi_c^{*,\pm}\|_{\mp \alpha}\leq D_{N}\nu^{N+1}. $$

This implies that we can extend the estimates to the whole half-line. Hence the approximate solution $\phi^{*,s,\pm}_{N}$ satisfies each of the required properties.

\end{proof}

Note that because we cannot obtain any good bounds for $D_N$, this approximate solution does not necessarily converge to an exact solution of the adjoint Orr-Sommerfeld equation.
\subsection{Approximate fast modes}\label{fast}
We will construct a fast growing and decaying approximate solution. These are linked, rather than to an Airy-type equation as for the original equation, to the equation
$\Delta_{\alpha}\l (\partial_{yy}-\mu_*^2)\phi\r=0. $ 
 We shall prove the following:

\begin{prop} Fix $N\in \mathbb{N}$. There exist two approximate solutions $\phi^{*,f,\pm}_N$ to the adjoint Orr-Sommerfeld equation \eqref{eq:orradj}, in the form 
\begin{equation}\label{eq:adjfast}\phi^{*,f,\pm}_N= e^{{\theta_{\pm}}/{\nu^{1/2}}}, \qquad \theta_{\pm}=\sum_{j=0}^{2N}\theta_{j,\pm}\nu^{j/2},\end{equation}
where in particular
\begin{equation} \theta_{0,\pm}=\pm \nu^{1/2}\int_0^{y}\mu_*(x)\,\mathrm{d}x, \qquad \mu_*=\sqrt{\alpha}\nu^{-1/2}\sqrt{\alpha \nu-i(U_s-\bar{c})}. \end{equation}
Then $\phi_N^{*,f,\pm}$ satisfies
$$|\orr_{c,\nu}^* \phi^{*,f,\pm}_N(y)|\leq D_N\nu^{N+1}|\phi^{*,f,\pm}_N(y)|,\qquad \forall y>0.$$
In particular, we have $\phi^{*,f,\pm}_N(0)=1$, and as $\nu\to 0$, $$\partial_y \phi^{*,f,\pm}_N(y)= \mu_*+O(1), \qquad \partial_{yy}\phi^{*,f,\pm}_N(y)= \mu_*^2+\mu_*'+O(1).$$
\end{prop}
Note that $\mu_*'(0)=0$ since $U_s'(0)=0$. Since $\mu_*=\bar{\mu}$, if we only consider $\theta_{0,\pm}$ then the fast solutions of the adjoint Orr-Sommerfeld equation have the same absolute value as the fast Orr-Sommerfeld solutions. The imaginary part has the opposite sign, meaning they are off-phase by an angle of $\pi$. However, ultimately the complex conjugate of an (approximate) Orr-Sommerfeld solution cannot be an (approximate) adjoint solution, as if $\orr_{c,\nu}\phi=\psi$ then
$$ \orr^*_{c,\nu}\bar{\phi}=\bar{\psi}+2U_s'\bar{\phi}'+U_s''\bar{\phi}\sim C\nu^{-1/2}.$$
Notice that this error becomes smaller as $y\to \infty$, and in general when the derivatives of $U_s$ are small.

Let us therefore construct two solutions of the above form for arbitrarily large $N>0$ by matching the terms with the same order in $\nu$. Differentiating the ansatz \eqref{eq:adjfast}, we obtain
\begin{align*}
\partial_{y}\phi^{*,f,\pm}_N&= \frac{\theta_{\pm}'}{\nu^{1/2}}\phi^{*,f,\pm}_N,\\
    \partial_{yy}\phi^{*,f,\pm}_N&= \l\frac{(\theta_{\pm}')^2}{\nu}+\frac{\theta_{\pm}''}{\nu^{1/2}} \r \phi^{*,f,\pm}_N,\\
    \nu\partial_y^4 \phi^{*,f,\pm}_N&=\l\frac{(\theta_{\pm}')^4}{\nu}+6 \frac{(\theta_{\pm}')^2\theta_{\pm}''}{\nu^{1/2}}+4\theta_{\pm}'\theta_{\pm}'''+3(\theta_{\pm}'')^2+\nu^{1/2}\theta_{\pm}'''' \r \phi^{*,f,\pm}_N.
\end{align*}
In the leading order equation (i.e. the terms of order higher than $\nu^{-1/2}$), the only difference is the sign change in the $\Delta_{\alpha}^2$ term and the conjugate on $c$, because the terms where at least one derivative hits $U_s-\bar{c}$ are of order $\nu^{-1/2}$ or lower. Recalling that we have assumed $|\alpha|\lesssim \nu^{-1/2}$, we thus obtain 
$$(U_s-\bar{c})\l(\theta_{0,\pm}')^2-\nu\alpha^2\r+\frac{1}{i\alpha}\l (\theta_{0,\pm}')^4-2\nu\alpha^2(\theta_{0,\pm}')^2+\nu^2\alpha^4\r=0.$$
Therefore,
\begin{equation}\label{eq:adjfast1} (\theta_{0,\pm}')^2=\nu\l \alpha^2 -i\alpha \frac{U_s-\bar{c}}{\nu} \r,\end{equation}
which has two solutions
$$\theta_{0,\pm}' = \pm \nu^{1/2}\mu_*, \qquad \mu_*=\sqrt{|\alpha|}\nu^{-1/2}\sqrt{\alpha \nu-i(U_s-\bar{c})},$$
where we recall that in $\mu_*$ we take the root with positive real part.
For the terms of order $\nu^{-1/2}$, some derivatives of $U_s-\bar{c}$ will also appear. More precisely, we have
$$4(\theta_{0,\pm}')^3\theta_{1,\pm}'-4\nu\alpha^2\theta_{0,\pm}'\theta_{1,\pm}'+2i\alpha(U_s-\bar{c})\theta_{0,\pm}'\theta_{1,\pm}'=-i\alpha U_s'\theta_{0,\pm}'+i\alpha (U_s-\bar{c})\theta_{0,\pm}''-6\theta_{0,\pm}'\theta_{0,\pm}''-2\alpha^2\nu\theta_{0,\pm}''$$
hence, by \eqref{eq:adjfast1}, 
$$\theta_{1,\pm}' = -\frac{-i\alpha U_s'\theta_{0,\pm}'+i\alpha (U_s-\bar{c})\theta_{0,\pm}''-6\theta_{0,\pm}'\theta_{0,\pm}''-2\alpha^2\nu\theta_{0,\pm}''}{2i\alpha (U_s-\bar{c}) }, $$
from which we conclude that $\theta_{1,\pm}'$ is bounded, uniformly in $\alpha =o(\nu^{-1/2})$ and $\nu \ll 1$. From here one constructs the other $\theta_{j,\pm}$ by induction as in \cite{grenieros}. The only difference is a sign change in the terms generated by ${\frac{\nu}{i\alpha}} \Delta_\alpha^2$ and some terms with derivatives of $U_s-\bar{c}$ (up to the second order)   appearing as known terms in each equation.

 Then, as in the proof of Proposition 2.1 in \cite{grenieros}, we conclude.

 \begin{remark}\label{fastder} Denote with $\phi^{f,\pm}_N$ the approximate fast growing or decaying solution of the original Orr-Sommerfeld equation constructed using the same procedure above. It follows from the above proof that if $U_s$ vanishes at the origin with its first $k$ derivatives, then $\partial_y^j \phi^{*,f,\pm}_N(0)=\partial^j \phi^{f,\pm}_N(0)$ for all $j\leq k-2$. In particular, if $k\geq 4$ then $\phi^{*,f,-}_N$ and $\phi^{f,-}_N$ satisfy the same Navier boundary condition.
 \end{remark}
 \begin{remark}
     The coefficients $\theta_{1,+}$ and $\theta_{1,-}$ satisfy
     $$\theta_{1,+}'+\theta_{1,-}'= \frac{6\theta_{0,\pm}'\theta_{j,\pm}''}{i\alpha(U_s-\bar{c})}=\frac{6\nu\mu_*\mu_*'}{i\alpha(U_s-\bar{c})}= \frac{-3U_s'}{U_s-\bar{c}}=O(e^{-\eta_0y}),\qquad y\to \infty.$$
     As a result, to leading order, as $\nu \to 0$ we have
     $$\phi^{*,f,+}\phi^{*,f,-}=O(e^{-\eta_0 y}).$$
 \end{remark}

\section{Eigenvalues and eigenmodes}\label{sec4}
In Definition \ref{rayeig} we introduced the notion of eigenvalues and eigenmodes for the Rayleigh equation. We also showed that the adjoint Rayleigh equation has the same eigenvalues, though the eigenmodes are different. In this section, we study the same concept for the viscous counterparts.  In Proposition \ref{iff}, we prove that the eigenvalues of the Orr-Sommerfeld equation and its adjoint must match, as in the inviscid case. We also show that for any Rayleigh eigenvalue $c_0$ there exists an integer $\kappa>0$ such that on a ball around $c_0$ there exist exactly $\kappa$ (adjoint) Orr-Sommerfeld eigenvalues $(c^j_{\nu},\nu)$  (counted with their multiplicity) with
$$|c_{\nu}^j-c_0|^{\kappa}\sim C_1\nu^{1/2}\left|\frac{1+C_2\alpha\nu^{\gamma}}{1+C_3\alpha\nu^{\gamma-1/2}}\right|,\qquad \nu \to 0^+,\quad C_1,C_2,C_3\in \mathbb{C},C_1\neq 0.$$
Note that, regardless of the value of $\gamma\in\mathbb{R}$, we always at least have $|c_{\nu}^j-c_0|^{\kappa}=O(\nu^{1/2})$. The rate of convergence improves as $\gamma \to -\infty$; the no-slip case corresponds to $\gamma \to +\infty$ for which the rate of convergence approaches $\nu^{1/2}$. Even for the original Orr-Sommerfeld equation, we could not find this result proven in the literature.
\subsection{Green function}\label{green}
For each triple $(\alpha,c,\nu)$, using the approximate fundamental solutions constructed in Section \ref{appsol}, we are looking to construct a Green function for \eqref{eq:orrs}, namely a function $G_{c,\nu}(x,y)\in C^2(\mathbb{R}_+\times \mathbb{R}_+)$ such that
\begin{equation}\label{eq:green1}
    \begin{cases}
    \textup{Orr}^*_{c,\nu}G_{c,\nu}(x,y)=\delta_{y}(x);\\ 
    G_{c,\nu}(x,0^+)=\lim_{y\to +\infty}G_{c,\nu}(x,y)=0;\\
    \partial_{yy}G_{c,\nu}(x,0^+)=\nu^{-\gamma}\partial_yG_{c,\nu}(x,0^+).
    \end{cases}
\end{equation}
More precisely, by the first line we mean that for any smooth and bounded function $\phi$, we have
$$\orr^*_{c,\nu}\int_0^{\infty}G_{c,\nu}(x,y)\phi(x)\,\mathrm{d}x=\phi(y).$$

We adopt a slightly different approach to \cite{grenieros}. We first construct an approximate Green function $G^I_{c,\nu,N}$ satisfying $\orr_{c,\nu}^* G^I_{c,\nu,N}=\delta_y(x)+\textup{err}_N(x,y)$, where $\textup{err}_N(x,y)$ is a function  of order $O(\nu^N)$ as $\nu\to 0$. We then correct this error to obtain an exact Green function $G_{c,\nu}^I$ solving $\orr^*_{c,\nu}G_{c,\nu}^I(x,y)=\delta_y(x)$ without the boundary condition \eqref{eq:green1}$_3$. This allows us to construct exact solutions for the adjoint Orr-Sommerfeld equation. Finally, we write
$$G_{c,\nu}= G_{c,\nu}^I+G_{c,\nu}^b, $$
where $G_{c,\nu}^b$ solves $\orr^*_{c,\nu}G_{c,\nu}^b=0$ and corrects the boundary conditions of $G^I_{c,\nu}$.  The correction of the boundary condition is only possible when $(c,\nu)$ is different from certain values, namely the eigenvalues of the $\orr^*_{c,\nu}$ operator.

Our construction of the Green function is inspired by the paper \cite{grenieros}. However, we offer an improvement to also include the case where $c$ is close or equal to a Rayleigh eigenvalue (as long as these differ from the adjoint Orr-Sommerfeld eigenvalues).

Define
\begin{equation}\label{eq:appgreen}
       G_{c,\nu,N}^I(x,y)= \begin{cases}a^{s,+}(x)\phi^{*,s,+}_N(y)+ a^{f,+}(x){\phi^{*,f,+}_N(y)} & y<x;\\ 
       a^{s,-}(x){\phi^{*,s,-}_N(y)}+a^{f,-}(x)\phi^{*,f,-}_N(y)& y>x;
    \end{cases}
\end{equation}
for some smooth functions $a^{s,\pm}, a^{f,\pm}$ to be derived. We introduce the $4\times 4$ matrix
$$M_{c,\nu}:=(M_{c,\nu}^i):=\begin{pmatrix}\partial_i \phi^{*,s,-}_N & \partial_i    \phi^{*,s,+}_N & \partial_i \phi^{*,f,-}_N & \partial_i \phi^{*,f,+}_N \end{pmatrix}, \qquad i=0,\dots,3. $$ 
We seek to prove that with an appropriate choice of coefficients $a^{s,\pm}$, $a^{f,\pm}$ we can apply Proposition \ref{generalgreen}. We thus need to construct $G_{c,\nu,N}^I(x,y)\in C^2(\mathbb{R}_+^2)$, smooth on $x\neq y$, and with a jump on the third derivative at $x=y$:
\begin{equation}\label{eq:jump} \partial_y^3 G_{c,\nu,N}^I(y_-,y)-\partial_y^3 G_{c,\nu,N}^I(y_+,y)=\frac{i\alpha}{\nu}.\end{equation}

In order to force these conditions we need to require
$$ (-a^{s,-},a^{s,+},-a^{f,-},a^{f,+})^T:= M_{c,\nu}^{-1}\begin{pmatrix} 0& 0 & 0 & -i\alpha/\nu\end{pmatrix}^T.$$
 By Proposition \ref{appsol}, 
in leading order as $\nu\to 0$ we have
\begin{equation}
M_{c,\nu}\sim \begin{pmatrix}\phi_{c}^{*,-} & \phi_c^{*,+} & \phi_{N}^{*,f,-} & \phi_N^{*,f,+} \\ 
\partial_y \phi_{c}^{*,-} & \partial_y\phi_{c}^{*,+} & -\mu_*\phi_N^{*,f,-} & \mu_*\phi_N^{*,f,+} \\ 
\partial_{y}^2\phi_{c}^{*,-} & \partial_{y}^2 \phi_{c}^{*,+} & \mu_*^2\phi_N^{*,f,-} & \mu_*^2\phi_N^{*,f,+} \\ 
\partial_{y}^3\phi_{c}^{*,-} & \partial_{y}^3 \phi_{c}^{*,+} & -\mu_*^3\phi_N^{*,f,-} & \mu_*^3 \phi_N^{*,f,+}
\end{pmatrix}.
\end{equation}
Note that the first two columns of the matrix are independent of the viscosity. 
We then have, as $\nu\to 0$,
$$\det M_{c,\nu}\sim 2\mu_*^5\phi_N^{*,f,-}\phi_N^{*,f,+}J_c^*,$$
where $J_c^*=J_c^*(y)$ is the determinant introduced in \eqref{eq:inviscid}. Recall that
$J_c^*$ is bounded away from zero.  Therefore, $M_{c,\nu}$ is invertible for $\nu$ small enough and
$$\begin{pmatrix} a^{s,-} \\ a^{s,+}\\ a^{f,-}\\ a^{f,+}\end{pmatrix} \sim \frac{-i\alpha}{\nu\mu_*^2}\begin{pmatrix}-\phi_c^{*,+}/J_c^*\\
-\phi_c^{*,-}/J_c^*\\
(2\phi_N^{*,f,-}\mu_*)^{-1}\\
(2\phi_N^{*,f,+}\mu_*)^{-1}
\\ \end{pmatrix}  $$
Also recall that the fast approximate solutions $\phi^{*,f,\pm}_N$ never vanish by construction. Since
$$\left|(\nu\mu_*^2)^{-1}\right|=|\alpha|^{-1}\left|\nu\alpha^2-i(U_s-\bar{c})\right|^{-1}\geq |\alpha|^{-1}|U_s-\bar{c}|^{-1}\geq |\alpha|^{-1}|\Im c|^{-1},$$ we conclude that $$\|a^{f,\pm}\|_{L^{\infty}}\sim \frac{C_{1,\pm}}{1+|\Im c|}(\mu_*)^{-1},\qquad \|a^{s,\pm}\|_{L^{\infty}}\sim \frac{C_{2,\pm}}{|\alpha|(1+|\Im c|)},\qquad \text{ as }\nu\to 0,$$  
for some $C_{1,\pm},C_{2,\pm}>0$ independent from $\alpha,c,\nu$. 
\begin{remark}\label{greencoeff}The choice of coefficients \eqref{eq:appgreen} is not unique, even for a fixed value of $N$. If we add a linear combination of $\phi^{*,s,{\pm}}_N$ and $\phi^{*,f,\pm}_N$ to $G^I_{c,\nu,N}$, it still satisfies all the required properties.
\end{remark}

\subsubsection{Exact Green function}\label{exactgreen}
In this section, we correct the error in the Green function constructed in the previous subsection. From this Section onwards, we assume we have fixed a specific choice of $N\in \mathbb{N}$, and never change it throughout.  As a result, we will eventually eliminate $N$ from the notation in order to simplify it. However, the construction that follows still depends on the choice of $N$ until a boundary condition for the Green function is fixed in Section \ref{greenboundary}.  

Let
$$\textup{err}_N(x,y):= \orr^*_{c,\nu}(G_{c,\nu,N}^I(x,y))-\delta_{x=y}(x,y).$$
Note that this is a function (rather than a distribution) in $x$ and $y$, and is of course smooth outside of the line $x=y$.
From the construction of $G_{c,\nu,N}^{I}$ we know there exists a constant $C_N>0$ such that
$$| \textup{err}_N(x,y)|\leq C_{N}\nu^N\left||\alpha|^{-1} e^{-|\alpha (x-y)|}+\mu_*^{-1}e^{-\int_x^y \mu_*(y)}\right|. $$
Hereafter, for any function $T=T(x,y)$ we use the notation $T\star$ to denote the operator acting on a function $f=f(x,y)$ as 
$$(T \star f)(x,y)=\int_0^{\infty} T(c,\nu)f(z,y)\,\mathrm{d}z. $$

 \begin{prop}\label{greenbounds}Let $U(y)$ be a shear profile satisfying \eqref{eq:eta}.
Let $N\in \mathbb{N}$, $c\in \mathbb{C}$. We can construct an exact Green function $G_{c,\nu}^I$ for the adjoint Orr-Sommerfeld equation such that for some $\theta_0\in (0,1)$ and a constant $C_N=C_N(\theta_0,N,c)$, continuous in $c$ as long as $\Im c \gg 0$, we have
 \begin{equation}\label{eq:gfest1}
     |\partial_y^jG_{c,\nu}^I(x,y)|\leq C_N\l |\alpha|^{j-1}e^{-\theta_0 |\alpha (x-y)|}+ |\mu_*(x)|^{j-1}e^{-\theta_0|\int_x^y \mu ^*|}\r+O\l {C_N\nu^N}\r,
 \end{equation}
 for $j=0,1,2,3$.
 \end{prop}
Note that $G_{c,\nu}^I$ may still depend on $N$, as we did not require it to satisfy any specific boundary condition.
\begin{proof}
Starting from $G_{c,\nu,N}^{I}$ we can construct the exact Green function as 
$$G_{c,\nu}^I(x,y): = \sum_{n=0}^{\infty}G_n(x,y),\qquad G_n:= (-\textup{err}_N \star)^n (G^{I}_{c,\nu,N}). $$
In other words, we have
$$G_n(x,y)=(-1)^n\int_0^{\infty}\textup{err}_N(c,\nu_n)\int_0^{\infty}\dots \int_0^{\infty}\textup{err}_N(z_2,z_1)G^I_{c,\nu,N}(z_1,y)\,\mathrm{d}z_1\dots \mathrm{d}z_n.$$

Let us prove that the series converges in $L^{\infty}$. Fix $j\in \lll 0;1;2;3\rrr$. Suppose by induction that $G_n$ satisfies
  \begin{equation}\label{eq:partial}|\partial_y^j G_n(x,y)|\leq  \l C_{N,\theta_0}\l |\alpha|^{j-1}+\|\mu_*\|_{\infty}^{j-1}\r\nu^{N}\r^n e^{-\theta_0 |\alpha (x-y)|}, \end{equation}
for some $\theta_0 \in (0,1)$. This is true for $n=0$ as long as $\nu$ small enough, which implies 
$$|\alpha(x-y)|\leq \int_x^y \Re \mu_*,\qquad \forall x, y\geq 0, $$
since $\Re \mu_*>0$, so we can absorb the fast decaying component. We then have
\begin{align*}
    |\partial_y^j G_{n+1}(x,y)|&\leq C_{N,\theta_0}^n \nu^{N}\l |\alpha|^{j-1}+\|\mu_*\|_{\infty}^{j-1}\r^{n+1} \int_0^{\infty}  e^{-\theta|\alpha (z-y)|} e^{-|\alpha (x-z)|} \,\mathrm{d}z\\
    &\leq C_{N,\theta_0}^n \nu^{N}\l |\alpha|^{j-1}+\|\mu_*\|_{\infty}^{j-1}\r^{n+1}e^{-\theta_0 |\alpha (x-y)|}C_{N,\theta_0},
\end{align*}
where
$$C_{N,\theta_0}=\sup_{x\geq 0} \int_0^{\infty}e^{-(1-\theta_0)|x-z|}\,\mathrm{d}z. $$
Since $G_0$ satisfies \eqref{eq:partial} with $\theta_0=1$ and in particular for any $\theta_0\in (0,1]$, we conclude by induction that \eqref{eq:partial} holds for all $n\in \mathbb{N}$, and thus \eqref{eq:gfest1} holds for any $\theta_0\in (0,1)$. 
So taking $\nu$ small enough, the series converges. Now, by construction we have $\orr^*_{c,\nu}G_{c,\nu,N}^I(x,y)=\delta_{x=y}(x,y)+\textup{err}_N(x,y)$. Therefore
 $$\orr_{c,\nu}\l \sum_{n=0}^M G_n\r = \sum_{n=0}^M (-1)^n (\textup{err}_N \star )^n (\delta_{x=y} + \textup{err}_N)=\delta_{x=y} + (-1)^M(\textup{err}_N  \star)^{M+1},$$
 which converges to $\delta_{x=y}$ as $M\to +\infty$. 
 \end{proof}
\subsection{Exact eigenvalues and eigenmodes}\label{exacteig}

It is well-known that in the no-slip case, shear flows which are unstable for the linearized Euler equation are also unstable for the linearized Navier-Stokes equations, as long as the viscosity is small enough. We extend this result to the adjoint Orr-Sommerfeld equation with the viscosity-dependent Navier boundary condition.

 Using the exact Green function, we may now construct the exact adjoint Orr-Sommerfeld solutions starting from the approximate ones. The exact solutions still depend on $N$, since we do not require any normalization, but we still drop the $N$ from the notation for simplicity.

\begin{prop}\label{exactsols}Let $N\in \mathbb{N}$, and let $\phi_N^{*,s,\pm}$, $\phi^{*,f,\pm}$,  be the four approximate solutions of the adjoint Orr-Sommerfeld equation constructed in Section \ref{sec3}. Then there exists four smooth independent exact solutions $\phi^{*,s,\pm}_{c,\nu}$, $\phi^{*,f,\pm}_{c,\nu}$ to the adjoint Orr-Sommerfeld equations, such that
\begin{equation}
    \phi^{*,k,\pm}_{c,\nu}\sim \phi^{*,k,\pm}_N(1-\tilde{D}_N\nu^{N+1}),\qquad \text{ as }\nu\to 0,
\end{equation}
where $\tilde{D}_N\neq 0$.
In particular, the solutions $\phi^{*,k,-}_{c,\nu}$ are decaying as $y\to \infty$, while the $\phi^{*,k,+}_{c,\nu}$ are growing.
\end{prop}
 
 \begin{proof}
 
 Fix $\alpha, c$ and $\nu$. Let $\phi^{*,\ell,-}_N$, $k=s,f$ be the approximate decaying and growing solutions. Define
\begin{equation}\label{eq:preimage} \tilde{\phi}^{k,-}_{N}(y):=\int_0^{\infty}G_{c,\nu}^I(x,y)\orr^*_{c,\nu}\phi^{*,k,-}_N(x)\,\mathrm{d}x,\end{equation}
and
$$\phi^{*,k,-}_{c,\nu}:=\phi^{*,k,-}_N-\tilde{\phi}^{k,-}_{N}.$$
We thus have
\begin{equation}\label{eq:exact}\phi^{*,k,-}_{c,\nu}\sim \phi^{*,k,-}_N\l 1- \tilde{D}_N\nu^{N+1}\r, \end{equation}
for some $\tilde{D}_N>0$ which is related to the constants from Proposition \ref{appsol} and Proposition \ref{greenbounds}. Then $\phi^{*,k,-}_{c,\nu}$ satisfies the following properties:
\begin{enumerate}
\item $\phi^{*,k,-}_{c,\nu}\neq 0$ at least for $\nu$ small enough (depending on $N$) by \eqref{eq:exact};
\item it has the same asymptotic behavior as $\nu \to 0$ as $\phi^{k,-}_N$, again by \eqref{eq:exact};
\item $\orr^*_{c,\nu}\phi^{*,k,-}_{c,\nu}=0$ because $\orr^*_{c,\nu}G_{c,\nu}^I(x,y)=\delta_y(x)$, by Proposition \ref{generalgreen};
\item it is smooth, because $\orr_{c,\nu}^* \phi_N^{*,k,-}$ is smooth and $G^I_{c,\nu}$ satisfies the assumptions of Proposition \ref{generalgreen}.
\end{enumerate}
This tells us that $\phi^{*,k,-}_{c,\nu}$ is a non trivial solution of the adjoint Orr-Sommerfeld equation, with the same asymptotic properties as the corresponding approximate solutions as $\nu\to 0$. Again, we point out that $\phi_{N}^{*,k,-}$ might not converge pointwise to $\phi^{*,k,-}_{c,\nu}$ as $N\to \infty$, due to the growth of $\tilde{D}_N$. 

The above procedure does not immediately apply to the growing solutions because the integral \eqref{eq:preimage} fails to converge. To treat this case, we require an adapted Green function. Define 
\begin{equation}\tilde{G}_{c,\nu,N}^I:=\begin{cases}0 & y<x;\\ 
   a^{s,-}(x){\phi^{*,s,-}_N(y)}+a^{f,-}(x)\phi^{*,f,-}_N(y)-a^{s,+}(x)\phi^{*,s,+}_N(y)- a^{f,+}(x){\phi^{*,f,+}_N(y)}& y>x.
    \end{cases}\end{equation}
By Remark \ref{greencoeff}, this is still an approximate Green function with the same coefficients as the ones derived in Section \ref{green}. It has a bounded support in $x$ for each fixed $y$, but in turn has fast growth as $y\to \infty$ for a fixed $x$. However, using the bounded support property we can still correct the error as in Section \ref{exactgreen}:
$$\tilde{G}_n(x,y):=(-\textup{err}_N\star)^n \tilde{G}_{c,\nu,N}^I(x,y)=(-1)^n\int_{0\leq z_n\leq \dots \leq z_1\leq y}\textup{err}_N(c,\nu_n)\dots \textup{err}_N(z_2,z_1)G(z_1,y)\,\mathrm{d}z_1\dots \mathrm{d}z_n, $$
so the integral in the right-hand side converges and is of order $(C_N\nu^N)^n$. In particular, from the above expression we also have $\tilde{G}_n(x,y)=0$ for $x>y$. Therefore for $\nu>0$ small enough the series $\sum_{n=0}^{\infty}\tilde{G}_n$ converges to an exact Green function $\tilde{G}_{c,\nu}^I(x,y)$ with $\tilde{G}_{c,\nu}^I(x,y)=0$ for all $x>y$. We can now set for $k=s,f$:
\begin{align*}\tilde{\phi}^{k,+}_{N}(y)&:=\int_0^{\infty}\tilde{G}_{c,\nu}^I(x,y)\orr_{c,\nu}^*\phi_N^{*,k,+}(x)\,\mathrm{d}x=\int_0^{y}\tilde{G}_{c,\nu}^I(x,y)\orr_{c,\nu}^*\phi_N^{*,k,+}(x)\,\mathrm{d}x;\\
\phi_{c,\nu}^{*,k,+}&:=\phi_N^{*,k,+}-\tilde{\phi}^{k,+}_{N}\sim \phi_N^{*,k,+}(1-\tilde{D}_N\nu^{N+1})\qquad \text{ as }\nu\to 0.
\end{align*}
Then $\phi_{c,\nu}^{*,k,+}$ is a smooth, exact solution to Orr-Sommerfeld which is growing as $y\to +\infty$, with at least the same order as $\phi_{N}^{*,k,+}$.
\end{proof}
Now that we have constructed the exact solutions, we can use them to extract more information on the viscous eigenvalues and eigenmodes. Recall the domain $\mc{D}_{\nu}$ defined in \eqref{eq:domain2}, which forces our Navier boundary condition.
\begin{deff}\label{orreigdef}
    Fix $\alpha\in\mathbb{Z}\setminus\lll 0\rrr$. A pair $(c,\nu)$ with $\Im c>0$ is called an \emph{adjoint Orr-Sommerfeld eigenvalue} if there exists $0\neq \phi \in \mc{D}_{\nu}$ with $\orr^*_{c,\nu}\phi=0$. In this case, $\phi$ is called an \emph{eigenmode}.
\end{deff}

Consider now the so-called \emph{Evans function} $\mc{E}_{c,\nu}$, defined as
\begin{equation}\label{eq:gc}E_{c,\nu}^{*}:=\begin{pmatrix}\phi^{*,s,-}_{c,\nu}(0) & \phi^{*,f,-}_{c,\nu}(0) \\ 
\partial_y\phi^{*,s,-}_{c,\nu}(0)-\nu^{\gamma}\partial_{yy}\phi^{*,s,-}_{c,\nu}(0) & \partial_y\phi^{*,f,-}_{c,\nu}(0)-\nu^{\gamma}\partial_{yy}\phi^{*,f,-}_{c,\nu}(0)
\end{pmatrix},\qquad \mc{E}^*(c,\nu):=\det E^{*}_{c,\nu}.\end{equation}
Then $(c,\nu)$ is an eigenvalue for the adjoint Orr-Sommerfeld equation if and only if $\mc{E}^{*}_{c,\nu}=0$. Indeed, by Proposition \ref{exactsols} there are only two independent decaying solutions, so $(c,\nu)$ is an eigenvalue if and only if there exist coefficients $a,b$ such that $\phi_{c_{\nu},\nu}^*=a\phi^{*,s,-}+b\phi^{*,f,-}$ is an eigenmode. In order to study the eigenvalues of the adjoint Orr-Sommerfeld equations, we thus study the zeros of $\mc{E}^*(c,\nu)$. The result that follows yields a differentiable curve of eigenvalues $\nu\mapsto (c_{\nu},\nu)$, for $\nu>0$ small enough. It is obtained through the implicit function theorem, which we must apply it to an adapted Evans function $\tilde{\mc{E}}^*_{c,\nu}$,  as $\mc{E}^*(c,\nu)$ is not even continuous as $\nu\to 0$. To this purpose, we define the function
$$\mc{O}_{\gamma}(\nu):=\lim_{c\to c_0}\frac{\partial_y\phi^{*,s,-}_{c,\nu}(0)-\nu^{\gamma}\partial_{yy}\phi^{*,s,-}_{c,\nu}(0)}{\partial_y\phi^{*,f,-}_{c,\nu}(0)-\nu^{\gamma}\partial_{yy}\phi^{*,f,-}_{c,\nu}(0)}.$$
Notice that the above quotient is continuous in $(c,\nu)$ around $(c_0,0)$. Furthermore,
$$\mc{O}_{\gamma}(\nu)\sim C_1\sqrt{\nu}\frac{1+C_2\alpha\nu^{\gamma}}{1+C_3\alpha\nu^{\gamma-1/2}}\; \text{ as } \nu\to 0^+,\qquad C_1,C_2,C_3\in \mathbb{C},C_1\neq 0.$$
\begin{prop}\label{evans}Let $c_0$ be a Rayleigh eigenvalue. Let $\kappa>0$ be the integer defined in Proposition \ref{newsol}). Then for each $\nu>0$ small enough there are exactly $\kappa$ adjoint Orr-Sommerfeld eigenvalues $(c_{\nu}^j,\nu)$ (counted with their multiplicity) in an appropriate neighborhood of $c_0$. In addition, the $c_{\nu}^j$ satisfy the asymptotics \begin{equation}\label{eq:eigas} |c_{\nu}^j-c_0|^{\kappa}\sim C|\mc{O}_{\gamma}(\nu)|,\qquad \nu \to 0^+,\quad j=1,\dots, \kappa,\end{equation}
for some constant $C>0$ independent from $\nu$.
\end{prop}
\begin{proof} We seek to apply Lemma \ref{implicit} to $\mc{E}^*(c,\nu)$. However, $\mc{E}^*(c,\nu)$ is not even continuous as $\nu\to 0^+$. Therefore, recalling that $\phi^{*,f,-}_{c,\nu}(0)=1$ by construction, we introduce
$$\tilde{\mc{E}}^*(c,\nu):=  \mc{E}^*(c,\nu)\l \partial_y \phi^{*,f,-}_{c,\nu}(0)- \nu^{\gamma}\partial_{yy}\phi^{*,f,-}_{c,\nu}(0)\r^{-1}=  \phi^{*,s,-}_{c,\nu}(0)- \frac{\partial_y\phi^{*,s,-}_{c,\nu}(0)-\nu^{\gamma}\partial_{yy}\phi^{*,s,-}_{c,\nu}(0)}{\partial_y\phi^{*,f,-}_{c,\nu}(0)-\nu^{\gamma}\partial_{yy}\phi^{*,f,-}_{c,\nu}(0)}. $$
Then
$$\tilde{\mc{E}}^*(c,\nu)=0 \iff \mc{E}^*(c,\nu)=0\iff (c,\nu) \text{ is an adjoint eigenvalue }.$$
By Proposition \ref{newsol}, the function $\tilde{\mc{E}}^*(c,\nu)$ is analytic in $\bar{c}$ and satisfies $\tilde{\mc{E}}(c_0,0)=0$. By the assumptions,
$$\phi_{c,\nu}^{*,s,-}(0)= \phi_c^-(0)+O(\nu)\sim A(\bar{c}-\bar{c}_0)^{\kappa} +O(\nu)\qquad A\in \mathbb{C},A\neq 0.$$
In conclusion, at $(c_0,0)$ we have an expansion of the type
$$ \tilde{\mc{E}}^*(c,\nu)\sim A(\bar{c}-\bar{c}_0)^{\kappa}+\mc{O}_{\gamma}(\nu)+o\l (\bar{c}-\bar{c}_0)^{\kappa}\r + o(\mc{O}_{\gamma}(\nu)),$$
By Lemma \ref{implicit} applied to $c\mapsto \tilde{E}^*(\bar{c},\nu)$ it then follows that for all $\nu>0$ small enough there exist $\kappa$ (counted with their multiplicity) zeroes $c_1(\nu),\dots c_{\kappa}(\nu)\to c_0$ as $\nu\to 0^+$, satisfying the asymptotics \eqref{eq:eigas} with $C=|A|^{-1}$.
\end{proof}

We conclude that near each Rayleigh eigenvalue $c_0$ there exists at least an adjoint Orr-Sommerfeld eigenvalue $c_{\nu}$ with $c_{\nu}-c_0=O(\nu^{1/{(2\kappa)}})$ as $\nu\to 0$. 
Given an adjoint eigenvalue $c_{\nu}$, we can write the associated adjoint Orr-Sommerfeld eigenmode as
$$\phi_{c_{\nu},\nu}^*\sim \phi^{*,s,-}_{c_{\nu},\nu}+C\mc{O}_{\gamma}(\nu)\phi^{*,f,-}_{c_{\nu},\nu}, $$
for some $C\in \mathbb{C}$ with $C\neq 0$. 
In conclusion, for any $\gamma\in \mathbb{R}$ the eigenfunction of the adjoint Orr-Sommerfeld operator is of the form (up to a multiplicative constant)
\begin{equation*} \phi^*_{c_{\nu},\nu}\sim \frac{\bar{\phi}^{-}_{c_{\nu},\nu}}{U_s-\bar{c}_{\nu}}+C\mc{O}_{\gamma}(\nu)\phi^{*,f,-}_{c_{\nu},\nu}=\frac{\bar{\phi}_{c_{\nu}}^-}{U_s-\bar{c}_{\nu}}+O(\mc{O}_{\gamma}(\nu))=\frac{\bar{\phi}_{c_0}^-}{U_s-\bar{c}_{0}}+O(\mc{O}_{\gamma}(\nu)),\end{equation*}
recalling that $\phi_{c_{\nu}}^-$ is the decaying fundamental solution of the Rayleigh equation.
 
\subsubsection{Boundary condition}\label{greenboundary}
In this Section, we correct the boundary condition in the Green function constructed in Section \ref{exactgreen}. Naturally, this introduces a singularity which corresponds to the adjoint Orr-Sommerfeld eigenvalues introduced in Definition \ref{orreigdef}.  Assuming $(c,\nu)$ is not such an eigenvalue, then the vectors $v_s:=\l\phi^{*,s,-}_{c,\nu}(0),\partial_y\phi^{*,s,-}_{c,\nu}(0)-\nu^{\gamma}\partial_{yy}\phi^{*,s,-}_{c,\nu}(0)\r$ and $v_f:=\l \phi^{*,f,-}_{c,\nu}(0),\partial_y\phi^{*,f,-}_{c,\nu}(0)-\nu^{\gamma}\partial_{yy}\phi^{*,f,-}_{c,\nu}(0)\r$ span $\mathbb{R}^2$. Hence, given $w(x):=\l G_{c,\nu}^I(x,0),\partial_y G_{c,\nu}^I(x,0)-\nu^{\gamma}\partial_{yy}G_{c,\nu}^I(x,0) \r $, we can always find a linear combination $b^s(x) {v}_s+b^f(x){v}_f$ such that $b^s(x) {v}_s+b^f(x){v}_f+w(x)=0$.
In fact, the vector $(b^s,b^f)$ is given by 
$$(b^s,b^f)^T =({E}_{c,\nu}^*)^{-1}(-w), $$
where
\begin{align*}(E_{c,\nu}^*)^{-1}&=(\mc{E}^*(c,\nu))^{-1}\begin{pmatrix}\partial_{y}\phi_{c,\nu}^{*,f,-}(0)-\nu^{\gamma}\partial_{yy}\phi_{c,\nu}^{*,f,-}(0) & -\phi_{c,\nu}^{*,f,-}(0)\\ 
-\partial_y \phi_{c,\nu}^{*,s,-}(0)+\nu^{\gamma}\partial_{yy}\phi_{c,\nu}^{*,s,-}(0) & \phi_{c,\nu}^{*,s,-}(0)\end{pmatrix}\\
&\sim  (\mc{E}^*(c,\nu))^{-1}\begin{pmatrix} \mu_*(C_1+C_2 \nu^{\gamma}\mu_*) & -1\\ C_3\alpha+C_4\alpha^2\nu^{\gamma} & C_5(\bar{c}-\bar{c}_0)^{\kappa}\end{pmatrix},\qquad (c,\nu)\to (c_0,0). 
\end{align*}
Therefore,
$$|(b^s,b^f)^T|\leq |\mc{E}^*(c,\nu)|^{-1}|w|. $$
This quantity will blow up with order $(c-c_{\nu})^{-\kappa}$ as $c$ approaches an adjoint Orr-Sommerfeld eigenmode $c_{\nu}$. \begin{comment}
In fact, we have
\begin{equation}
    \begin{pmatrix}b^s \\ b^f \end{pmatrix} \sim \frac{1-\nu^{\gamma-1/2}\mu}{(c-c_{\nu})\nu^{1/2}\mu^2}\begin{pmatrix}-1 \\ \phi_c^-(0)\end{pmatrix}.
\end{equation}
\end{comment}
We can then take
$$G_{c,\nu}^b(x,y):= b^s(x)\phi^{*,s,-}_{c,\nu}(y)+b^f(x)\phi^{*,f,-}_{c,\nu}(y),\qquad G_{c,\nu}(x,y):=G^I_{c,\nu}(x,y)+G^b_{c,\nu}(x,y).$$
Since $\orr_{c,\nu}G_{c,\nu}^b=0$, the equation will still be satisfied. The coefficients $b^s,b^f$ are continuous with respect to $c$ as long as we are away from any adjoint Orr-Sommerfeld eigenvalue. We can therefore conclude the following.

 \begin{theorem}\label{gfbounds}
 Let $(c,\nu)$ be such that $\mc{E}^*(c,\nu)\neq 0$, i.e. $(c,\nu)$ is not an adjoint Orr-Sommerfeld eigenvalue. Let $G_{c,\nu}(x,y)$ be the exact Green function for the adjoint Orr-Sommerfeld equation satisfying \eqref{eq:green1}. Then there exists $\theta_0\in (0,1)$ and a constant $C_N=C_N(\theta_0,c)$, continuous in $c$ as long as $\Im c \gg 0$, such that for all $\alpha\leq \nu^{-\zeta}$, $\zeta <\frac{1}{2}$ we have
 \begin{equation}\label{eq:gfest}
     \frac{|\partial_y^jG_{c,\nu}(x,y)|}{1+|\mc{E}^*(c,\nu)|^{-1}}\leq C_N\l |\alpha|^{j-1}e^{-\theta_0 |\alpha (x-y)|}+ |\mu_*(x)|^{j-1}e^{-\theta_0|\int_x^y \mu _*|}\r+O\l {C_N\nu^N}\r,\quad j=0,1,2,3.
 \end{equation}
 \end{theorem}
In conclusion, starting from different values of $N$, we have obtained different approximate Green functions $G_{c,\nu,N}$ which we could use to construct an exact Green function $G_{c,\nu}$ as above. This Green function, unlike the others, must be independent from $N$ because, as long as $(c,\nu)$ is not an exact adjoint Orr-Sommerfeld eigenvalue, the operator $\orr_{c,\nu}^*$ is injective by definition.
\subsection{Image and kernel of the Orr-Sommerfeld operator}
Thanks to our analysis of the adjoint operator, it is possible to obtain some useful information on the image of the original Orr-Sommerfeld operator. Throughout this section, let $\lll (c_{\nu},\nu)\rrr_{\nu>0}$ and $\lll (c^*_{\nu},\nu)\rrr_{\nu>0}$ be a family of eigenvalues of Orr-Sommerfeld (respectively, its adjoint) converging to a Rayleigh eigenvalue $c_0$ with $\Im c_0>0$ at a rate of $O(\sqrt{\nu})$, as predicted by Proposition \ref{evans}.

We know that  $\ker \orr^*_{c,\nu}\cap \mc{D}_{\nu}\subset\Img \orr_{c,\nu}^{\perp} $.
However, $\ker \orr^*_{c,\nu}$ is either trivial or is spanned by the eigenmode $\phi^*_{c,\nu}$ constructed in Section \ref{exacteig} when $c$ is an adjoint eigenvalue. But $\phi^*_{c,\nu}\in \mc{D}_{\nu}$ as it is smooth and satisfies the boundary conditions, so we conclude
$$ \ker \orr^*_{c,\nu}\subset\Img \orr_{c,\nu}^{\perp},\qquad (\ker \orr_{c,\nu}^*)^{\perp}\supset \overline{\Img \orr_{c,\nu}},$$
and the same where the adjoints are exchanged between the two sides of the inclusions. In particular, if $(c_{\nu},\nu)$ and $(c^*_{\nu},\nu)$ are eigenvalues and adjoint eigenvalues respectively, we have $\phi_{c_{\nu},\nu}\in\l \Img \orr^*_{c_{\nu},\nu}\r ^{\perp}$ and $\phi^*_{c^*_{\nu},\nu}\in \Img \orr_{c_{\nu}^*.{\nu}}^{\perp}$. Thanks to the Green function constructed in Section \ref{green}, we in fact deduce that the adjoint eigenvalues are precisely the eigenvalues of the original operator. Thus this property carries over from the Rayleigh operator and its adjoint, even though in the viscous case we could not find any explicit formula to derive an adjoint eigenmode from the original eigenmode.
\begin{prop}\label{iff}Fix $\alpha\in\mathbb{Z}$. A pair $(c,\nu)$ is an Orr-Sommerfeld eigenvalue if and only if it is an adjoint eigenvalue. 
\end{prop}
\begin{proof} We prove one direction, the other is identical. Suppose $(c,\nu)$ is an Orr-Sommerfeld eigenvalue, but not an adjoint eigenvalue. Then we can construct the adjoint Green function with boundary conditions $G_{c,\nu}$ as per Section \ref{greenboundary}. Let $\psi\in\mc{D}_{\nu}$ be the eigenmode for $\orr_{c,\nu}$, and let $\phi(y):=\int_0^{\infty}G_{c,\nu}(x,y)\psi(x)\,\mathrm{d}x$. As $G_{c,\nu}$ satisfies the boundary conditions, so does $\phi$, and by Proposition \ref{generalgreen} we have $\orr^*_{c,\nu}\phi=\psi$. Since both $\phi$ and $\psi$ belong to $\mc{D}_{\nu}$, we then have
$$0< \lpl \psi,\psi\rer_{L^2} = \lpl \phi, \orr_{c,\nu}\psi\rer_{L^2} = 0, $$
a contradiction. 
\end{proof}

Exact Orr-Sommerfeld eigenvalues can be directly constructed in the same way we did for adjoint eigenvalues in Proposition \ref{evans}. But to first order approximation as $\nu\to 0$ and $c\to c_0$, the matrix $\mc{E}^*(c,\nu)$
corresponding to the Orr-Sommerfeld operator would simply be the complex conjugate of $E^{*}_{c,\nu}$. Thus the same would be true for its determinant, and hence its zeros i.e. the eigenvalues would remain unchanged. However, this is not a complete proof as in principle the lower order terms might change the picture. Regrettably, we were unable to find a direct and rigorous proof of Proposition \ref{iff}.

\begin{remark}
    From now on, we will simply refer to (adjoint) Orr-Sommerfeld eigenvalues as \emph{eigenvalues} without distinguishing the adjoint and original operator.
\end{remark}
The following result shows that eigenmodes (resp. adjoint eigenmodes) also do not belong to the image of the Orr-Sommerfeld operator (resp. adjoint) in correspondence of the eigenvalues. This can be seen as the viscous counterpart of Proposition \ref{rayim}.

\begin{prop}\label{image1} For $\nu>0$ small enough, we have $\bar{\phi}_{c_{\nu},\nu} \notin \Img \orr_{c_{\nu},\nu}$ and $\bar{\phi}^*_{c_{\nu},\nu}\notin \Img \orr^*_{c_{\nu},\nu}$. Furthermore, if $\phi_{c_0}^-\notin \Img \ray_{c_0}$, then $\phi_{c_{\nu},\nu}\notin \Img \orr_{c_{\nu},\nu}$ and $\phi_{c_{\nu},\nu}^*\notin \Img \orr^*_{c_{\nu},\nu}$ for $\nu$ small enough.
\end{prop}
\begin{proof}
Let us prove that $\phi_{c_{\nu},\nu}\notin (\ker \orr^*_{c_{\nu},\nu})^{\perp}$. We know that $c_{\nu}= c_0+O(\nu^{1/2})$ is an eigenvalue of $\orr_{c,\nu}$ as well as $\orr^*_{c,\nu}$. Let $\phi_{c_{\nu},\nu}^*\in \ker \orr^*_{c,\nu}$ be the corresponding eigenmode. Then $\phi_{c_{\nu},\nu}^*$ satisfies the boundary conditions by construction, and 
\begin{equation}\label{eq:orradj2}\phi^*_{c_{\nu},\nu}= \frac{\overline{\phi_{c_{\nu},\nu}}}{U_s-\overline{c}_{\nu}}+O(\nu^{1/2}), \end{equation}
with $\orr^*_{c_{\nu},\nu}\phi^*_{c_{\nu},\nu}=0$. Moreover, $\phi_{c_{\nu},\nu}^*$ spans $\ker \orr_{c_{\nu},\nu}^*$. We then have, for some constant $C>0$ independent from $\nu$,
\begin{align*}
     \left|\int_0^{\infty}\phi_{c_{\nu},\nu}{\phi}^*_{c_{\nu},\nu}\right|\geq Im c_{\nu}\int_0^{\infty}\frac{|\phi_{c_{\nu},\nu}|^2}{|U_s-c_{\nu}|^2}-C\nu^{1/2}\int_0^{\infty}|\phi_{c_{\nu},\nu}|,
\end{align*}
which is nonzero for $\nu$ small enough, since $\|\phi_{c_{\nu},\nu}\|_{L^1}$ is bounded uniformly in $\nu$. Hence $\bar{\phi}_{c_{\nu},\nu}\notin \l\ker {\orr^*_{c_{\nu},\nu}}\r^{\perp}\supset \overline{ \Img \orr_{c_{\nu},\nu}}$. The same inequality also tells us that $\bar{\phi}^*_{c_{\nu},\nu}\notin \ker \orr_{c_{\nu},\nu}^{\perp}$, and hence $\bar{\phi}^*_{c_{\nu},\nu}\notin \Img \orr^*_{c_{\nu},\nu}$.

As for the eigenmodes themselves, we know from Remark \ref{eigenimage} that 
$$ \phi_{c_0}^-\notin \Img \ray_{c_0}\implies \int_0^{\infty}\frac{(\phi_{c_0}^-)^2}{U_s-c_0}\neq 0.$$
In this case, using $\phi_{c_{\nu},\nu}=\phi_{c_0}^-+ O(\sqrt{\nu})$ and $(c_{\nu}-c_0)^{\kappa}=O(\sqrt{\nu})$ we conclude that for $\nu$ small enough
$$ \int_0^{\infty} \frac{(\phi_{c_{\nu},\nu})^2}{U_s-c_{\nu}}\neq 0 \implies \phi_{c_{\nu},\nu}\notin \Img \orr_{c_{\nu},\nu}.$$
For $\phi_{c_{\nu},\nu}^*$ the proof is the same.
\end{proof}
We have thus seen that when $(c,\nu)$ is an eigenvalue the Green function cannot exist, as the existence of a Green function strictly implies that $\mc{D}_{\nu}$ must be contained in the codomain.

However, we can always construct a Green function which only corrects one of the two boundary conditions. That allows us to characterize the image of $\orr^*_{c_{\nu},\nu}$ as the set of functions which are orthogonal to the Navier boundary condition of the Green function with respect to the $L^2$ product.

More precisely, recalling that $\phi_{c,\nu}^{*,f,-}(0)\neq 0$ for all $(c,\nu)$, we can define
\begin{equation}
    \tilde{G}^b_{c,\nu}(x,y):= -\frac{G_{c,\nu}^I(x,0)}{\phi^{*,f,-}_{c,\nu}(0)}\phi^{*,f,-}_{c,\nu}(y),
\end{equation}
so that $\tilde{G}_{c,\nu}:=G^I_{c,\nu}+\tilde{G}^b_{c,\nu}$ is a Green function for the adjoint Orr-Sommerfeld operator, satisfying the assumptions of Proposition \ref{generalgreen} and the boundary condition $G_{c,\nu}(x,0)=0$ for all $x>0$. 
\begin{prop}Let $(c,\nu)$ be an Orr-Sommerfeld eigenvalue, and let $\psi \in C^{3}\cap L^2(\mathbb{R}_+)$. Suppose that
\begin{equation}\label{eq:orthogorr}\int_0^{\infty}(\partial_y-\nu^{\gamma}\partial_{yy})\tilde{G}_{c,\nu}(x,0)\psi(x)\,\mathrm{d}x=0.
\end{equation}
Then $\psi \in \Img \orr_{c,\nu}^*$. In particular,  $\Img \orr^*_{c,\nu}$ has codimension one.
\end{prop}
\begin{proof}Defining
$$\phi(y):=\int_0^{\infty}\tilde{G}_{c,\nu}(x,y)\psi(x)\,\mathrm{d}x,$$
then by Proposition \ref{generalgreen} we have $\orr^*_{c,\nu}\phi=\psi$. Moreover, $\phi(0)=0$ since $\tilde{G}_{c,\nu}(x,0)=0$ for all $x$, and $\phi'(0)-\nu^{\gamma}\phi''(0)=0$ by \eqref{eq:orthogorr}. Thus $\phi \in \mc{D}_\nu$ and $\psi \in \Img \orr_{c,\nu}^*$. This implies that
$$\lpl \overline{(\partial_y-\nu^{\gamma}\partial_{yy})\tilde{G}_{c,\nu}(x,0)}\rer^{\perp}\cap C^{3}(\mathbb{R}_+)\subset \Img \orr^*_{c,\nu},$$
and hence $\Img \orr^*_{c,\nu}$ has at most codimension one, since $C^{3}\cap L^2$ is dense in $L^2$. But we already know from Proposition \ref{image1} that the codimension is at least one, so we conclude.
\end{proof}
Of course, we could have alternatively first corrected the Navier boundary condition. The resulting Green function at $y=0$ would then be equal to $(\partial_y-\nu^{\gamma}\partial_{yy})\tilde{G}_{c,\nu}(x,0)$ up to a multiplicative constant independent from $x$. Thus its orthogonal would define the same space, avoiding any contradiction.
\appendix

\section{Estimates on Green functions}
This Proposition ensures the validity of each of the Green functions considered in this paper. 
\begin{prop}\label{generalgreen}
    Let $n\geq 2$. Consider a linear differential operator of order $n$ given by
    $$ T:=\sum_{k=0}^ng_k(y) \pdv{^k}{y^k},$$
    where $g_k\in C^{\infty}([0,\infty);\mathbb{C})$ are smooth, bounded with all their derivatives as well as
    $\inf_{[0,\infty)}|1/g_n|>0.$ Suppose that there exists a function $G:\mathbb{R}_+\times \mathbb{R}_+\to \mathbb{C}$ given by  $$G(x,y)=\begin{cases} G_1(x,y) & y\geq x \\
    G_2(x,y) & x<y \end{cases}
  $$
  such that $G \in C^{n-2}([0,\infty)\times [0,\infty))$, while $G_1,G_2$ are smooth up to the boundary $x=y$, satisfy $T G_1(x,y)=T G_2(x,y)=0$ for $x\neq y$ and
    \begin{equation}\label{eq:cond}
           \partial_2^{n-1}G(y_-,y)-\partial_2^{n-1}G(y_+,y)=\frac{1}{g_n(y)},\qquad \forall y \in [0,\infty).
    \end{equation}
    Then for all $\phi\in C^{n-1}([0,\infty)];\mathbb{C})$ such that $\int_0^{\infty}|\partial_y^j G(x,y)\phi(x)|\,\mathrm{d}x<\infty$ for all $j=0,\dots, n-1$,
     the following properties hold:
    \begin{enumerate}
        \item We have $TG(x,y)=\delta_y(x)$, in the sense that for all $\phi$ in the above class we have
        \begin{equation}\label{eq:general1}
       T\int_0^{\infty} G(x,y)\phi(x)\,\mathrm{d}x=\phi(y).
\end{equation}
        \item Let $\phi$ in the above class, with $\phi \in C^{\max\lll n-1;k-n\rrr}([0,\infty);\mathbb{C})$, and let $f(y):=\int_0^{\infty} G(x,y)\phi(x)\,\mathrm{d}x$. Then  for all $y\in [0,\infty)$ and $k\in \mathbb{N}$ we have
    \begin{equation}\label{eq:general2}
        |f^{(k)}(y)|\leq C_k\left[ \max_{0\leq j\leq k-n}|\phi^{(j)}(y)|+\max_{0\leq j\leq n-1\land k}\left|\int_0^{\infty}  \partial_y^j G(x,y)\phi(x)\,\mathrm{d}x\right|\right],
    \end{equation}
    for some $C_k$ depending on the supremums over $\mathbb{R}_+$ of the first $k-n\lor 0$ derivatives of the coefficients $g_1,\dots,g_{n-1}$ and $1/g_n$.
    \end{enumerate}
\end{prop}
\begin{proof}
Split the integral as
$$T\int_0^{\infty} G(x,y)\phi(x)\,\mathrm{d}x=I_1+I_2+I_3,$$
where
$$ I_1=T\int_{y+\eps}^{\infty}G_1(x,y)\phi(x)\,\mathrm{d}x,\qquad I_2=T\int_{0}^{y-\eps}G_2(x,y)\phi(x)\,\mathrm{d}x,\qquad  I_3 = T\int_{y-\eps}^{y+\eps}G(x,y)\phi(x)\,\mathrm{d}x$$
For $I_3$, we know that $y\mapsto \int_{y-\eps}^{y+\eps}G(x,y)\phi(x)\,\mathrm{d}x$ is of class $C^{n-2}$ with a piecewise continuous $(n-1)$-th derivative. Therefore, whenever we apply an operator of order $\leq n-1$ to this term, the result will vanish as $\eps \to 0$. As for the term of order $n$, we have
\begin{equation}\label{eq:n1} \partial_y^{n}\int_{y-\eps}^{y+\eps}G(x,y)\phi(x)\,\mathrm{d}x= \partial_y \int_{y-\eps}^{y+\eps}\partial_y^{n-1}G(x,y)\phi(x)\,\mathrm{d}x+\sum_{j=1}^{n-1}\dv{^j}{y^j}\left[\partial_y^{n-j-1}G(x,y)\phi(x)\right]_{x=y-\eps}^{x=y+\eps}.\end{equation}
Let us focus on the first term in the right hand side above. Note that we cannot apply Leibniz's rule to differentiate, because the derivative of $\partial_y^{n-1}G$ is not a function. However, by \eqref{eq:cond} and since $G_1,G_2$ are smooth up to the boundary we can write
$$\partial_y^{n-1}G(x,y)=\tilde{G}(x,y)+ \frac{1}{g_n(y)}\chi_{\lll y>x\rrr},$$
where $\tilde{G}$ is continuous and has a locally bounded, and hence locally integrable, derivative. Therefore, as $\eps \to 0$ we have
$$ \lim_{\eps \to 0}\partial_y \int_{y-\eps}^{y+\eps}\partial_y^{n-1}G(x,y)\phi(x)\,\mathrm{d}x =\lim_{\eps \to 0} g_n(y)\partial_y \frac{1}{g_n(y)}\int_{y-\eps}^y 1\,\mathrm{d}x=0.$$

It remains to consider $I_1$ and $I_2$. We know that $TG_1(x,y)=TG_2(x,y)=0$ within the respective domains, therefore after applying $T$ the only surviving terms are those where one derivative hits the integration domain, according to Leibniz's rule. These terms are of order between $0$ and $n-1$ in $G(x,y)$. However, all the terms of order $\leq n-2$ are continuous in $y$ and thus vanish in the sum $I_2+I_3$ as $\eps \to 0$. It remains to consider the terms of order $n-1$ in $G(x,y)$. These must come from the highest order term $g_n\pdv{^n}{y^n}$, so they are of the form
$$-g_n(y)\sum_{j=1}^{n-1}\dv{^j}{y^j} \left[\partial_y^{n-j-1}G(x,y)\phi(x)\right]_{x=y-\eps}^{x=y+\eps}-g_n(y)\left[\partial_y^{n-1}G(x,y)\phi(x)\right]_{x=y-\eps}^{x=y+\eps}.$$
The terms in the sum will cancel out with the corresponding terms from $I_1$ as in \eqref{eq:n1}, so
$$I_2+I_3 \to g_n(y)\l -\partial_2^{n-1} G(y_+,y)+\partial_2^{n-1}G(y_-,y)\r \phi(y)=\phi(y),$$
and we conclude \eqref{eq:general1}.

To obtain \eqref{eq:general2}, we differentiate $k$ times the equality 
$$f(y)=\int_0^{\infty}G(x,y)\phi(x)\,\mathrm{d}x.$$
 For the first $n-1$ derivatives, all the derivatives get through the integral to $G$, so the estimate follows, without requiring any derivative of $\phi$. For higher order derivatives, by \eqref{eq:general1} we know that
$$ f^{(n)}(y)=\frac{1}{g_n(y)}\l \phi(y)-\sum_{j=0}^{n-1}g_j(y)\partial_y^jf(y)\r.$$
By applying this rule repeatedly, we ultimately only get derivatives of $\phi$ at $y$ of order up to $k-n\lor 0$, as well as derivatives of $f$ of order up to $n-1$ which are in the previous form. We conclude \eqref{eq:general2}.
\end{proof} 
\section{Treatment of the original operators}\label{original}
In this Appendix, we explain how the techniques of this paper can be adapted to the original Rayleigh and Orr-Sommerfeld operator. The adjoint Rayleigh operator
$$\ray^*_{c}= (U_s-\bar{c})\Delta_{\alpha}+2U_s'\partial_y$$
only differs from 
$$\ray_c=(U_s-c)\Delta_{\alpha}-U_s''\phi$$ by the complex conjugate on $c$ and the lower order term which is $2U_s'\phi'$ rather than $-U_s''\phi$.
To replicate the results of Section \ref{sec2}, we first consider the operator
\begin{equation*}T_{\alpha;y_0}\phi(y)=\frac{1}{2\alpha}\int_{y_0}^{\infty}e^{-\alpha|y-x|}\frac{U_s''(x)}{U_s(x)-c}\phi(x)\,\mathrm{d}x.
\end{equation*}
With this definition, Lemma \ref{Top} carries over, where in \eqref{eq:yous} the point $y_k>0$ depends on the first $k$ derivatives of $U_s''/(U_s-c)$. From this, Proposition \ref{initial} and \ref{imageray}  follow. In particular, we have two smooth fundamental solutions
$$\phi_c^{\pm}(y) \sim e^{\pm \alpha y},\qquad y\to \infty.$$

To obtain Proposition \ref{newsol}, we similarly modify the $S_{c;c_0}$ operator to
$$S_{c;c_0}(\psi)=(c-c_0)\int_0^{\infty}G_{c_0}(x,y)\frac{U_s''(x)}{U_s(x)-c}\psi(x)\,\mathrm{d}x\implies \ray_{c_0} S_{c;c_0}\psi=(c-c_0)\frac{U_s''}{U_s-c}\psi.$$
In this case, we do not need to apply integration by parts as $\psi$ already appears in order $0$ in the right hand side. Estimate \eqref{eq:sest} then follows from Proposition \ref{generalgreen} and so does Proposition \ref{newsol}. With the exact solutions of the Rayleigh equation constructed, we can now construct the corrisponding Green functions and deduce Proposition \ref{raybounds} and \ref{ivpest} using Proposition \ref{generalgreen}.\\

In the viscous case, when comparing
$$\orr^*_{c,\nu}=(U_s-\bar{c})\Delta_{\alpha}+2U_s'\partial_y+\frac{\nu}{i\alpha}\Delta_{\alpha}^2$$
and 
$$\orr_{c,\nu}=(U_s-c)\Delta_{\alpha}-U_s'' -\frac{\nu}{i\alpha}\Delta_{\alpha}^2$$
we notice that, on top of the differences in the corresponding terms originating from the Rayleigh equation, the viscous term also has the opposite sign. The result is that for the approximate fast solutions, we have as $\nu\to 0$
$$\phi^{f,\pm}_N(y)\sim e^{\pm \int_0^y \mu_f}, \qquad \mu= |\alpha|^{1/2}\nu^{-1/2}\sqrt{\alpha\nu+i(U_s-c)} ,$$
where we notice that in $\mu$ the term $i(U_s-c)$ has the opposite sign compared to $\mu_*$, and of course also loses the conjugate on $c$. In first order approximation as $\nu\to 0$, the approximate fast solutions are therefore given by the complex conjugates of the corresponding adjoint solutions, and the growth or decay as $y\to \infty$ remains unchanged. As for the slow solutions, they are still readily constructed by approximation from the Rayleigh solutions and thus exhibit the same asymptotic behavior. 

Once all the approximate solutions are constructed, the results from Section \ref{sec4} consequently follow. In particular, the Evans function $\mc{E}(c,\nu)$ corresponding to \eqref{eq:gc} will be the complex conjugate of $\mc{E}^*(c,\nu)$
in first order approximation as $\nu\to 0$, and hence we find that its zeroes, i.e. the eigenvalues, satisfy the same asymptotics as the one obtained in Proposition \ref{evans}. This is in agreement with Proposition \ref{iff}, which states that eigenvalues and adjoint eigenvalues are in fact the same.

\section{Perturbation by a parameter of zeroes of complex analytic functions}
This technical Lemma is required in the proof of the existence and asymptotic properties of the Orr-Sommerfeld eigenvalues, Proposition \ref{evans}.
\begin{lemma}\label{implicit} Let $ \mathcal{U}\times[0,\nu_0]\subset \mathbb{C}\times \mathbb{R}$ a one-sided in $\nu$ neighborhood of $(c_0,0)$, and $f=f(c,\nu):\mathcal{U}\times[0,\nu_0]\to \mathbb{C}$ holomorphic in the first  variable and satisfying an expansion of the following type near $(c_0,0)$:
    $$f(c,\nu)= A (\bar{c}-\bar{c}_0)^{\kappa} +g(\nu)+R(c,\nu),$$
    where $$A\in \mathbb{C}\setminus \lll 0\rrr,\qquad |R(c,\nu)|=o((\bar{c}-\bar{c}_0)^{\kappa})+o(g(\nu))\quad \text{as }(c,\nu)\to (c_0,0),$$
    and $g:[0,\nu_0]\to \mathbb{C}$ satisfies $g(0)=0$ but $g(\nu)\neq 0$ for $\nu\neq 0$, and is continuous at $\nu=0$.

Then for all $\nu\in (0,\nu_0]$ small enough there exist exactly $\kappa$ zeroes (counted with their multiplicity) $c_1(\nu),\dots,c_{\kappa}(\nu)$ for $f$ in an appropriate neighborhood of $c=c_0$, with 
\begin{equation}\label{eq:zeroest}.
 |c_j(\nu)-c_0|^{\kappa}\sim  \frac{1}{|A|}|g(\nu)|,\qquad j=1,\dots,\kappa,\;\text{ as } \nu\to 0^+.
\end{equation}

\end{lemma}
\begin{proof} Fix $\eps>0$, and let $0<\alpha < \eps/(2+\eps)$.
   By reducing the value of $\nu_0$ if necessary, there exists $r_0>0$ such that for all $\nu \in [0,\nu_0],|c-c_0|\leq r_0$ we have
    $$|R(c,\nu)|\leq  \alpha\l |A||c-c_0|^{\kappa}+|g(\nu)|\r. $$
    Now fix any $\nu\in (0,\nu_0]$ small enough so that
    $$R:=\l \frac{1+\eps}{|A|}|g(\nu)|\r^{1/\kappa}\leq r_0.$$
    Then for all $|c-c_0|=R$ we have
    \begin{align*}|R(c,\nu)|&\leq \alpha\l |A||c-c_0|^{\kappa} +|g(\nu)|\r=\alpha |A|\l 1+\frac{1}{1+\eps}\r R^{\kappa}<|A|\l 1-\frac{1}{1+\eps} \r R^k\\
&=|A||c-c_0|^{\kappa}-|g(\nu)|\leq |A(\bar{c}-\bar{c}_0)^{\kappa}+g(\nu)|.\end{align*}
    Next, let 
    $$r:= \l \frac{1-\eps}{|A|}|g(\nu)|\r^{1/\kappa}.$$
    Then for $|c-c_0|=r$ we have
\begin{align*} |R(c,\nu)|\leq \alpha \l |A||c-c_0|^{\kappa}+|g(\nu)|\r= \alpha\l 2-\eps \r |g(\nu)|<\eps |g(\nu)|=|g(\nu)|-|A||c-c_0|^{\kappa}\leq |A(\bar{c}-\bar{c}_0)^{\kappa}+g(\nu)|.\end{align*}
    By Rouché's theorem from complex analysis, we deduce that $c\mapsto f(c,\nu)$ has the same number of zeroes of $c\mapsto A(\bar{c}-\bar{c}_0)^{\kappa}+g(\nu)$ in the annulus $\lll c:r<|c-c_0|<R\rrr$. Hence $f$ has exactly $\kappa$ zeroes $(c_1(\nu),\nu),\dots,(c_{\kappa}(\nu),\nu)$ in the annulus, counted with their multiplicity. Furthermore, we can also apply Rouché's theorem on $B_R(c_0)$, so these must be the only zeroes in $B_R(c_0)$ for each $\nu$. Using the explicit expressions for $r$ and $R$, we deduce
$$ \frac{1-\eps}{|A|}|g(\nu)|\leq |c_j(\nu)-c_0|^{\kappa}\leq \frac{1+\eps}{|A|}|g(\nu)|,\qquad j=1,\dots,\kappa,$$
for $\nu<\nu_0$ (depending on $\eps$). Since $\eps>0$ is arbitrary, we obtain \eqref{eq:zeroest}.
\end{proof}
 
\end{document}